\documentclass[aip,jcp,reprint]{revtex4-1}
\usepackage{dcolumn}
\usepackage{bm}
\usepackage[T1]{fontenc}  
\usepackage{cleveref}
\usepackage{amsmath}
\usepackage{amssymb}
\usepackage{booktabs}
\usepackage{lmodern}
\usepackage{standalone}
\RequirePackage{url}
\RequirePackage{graphicx}
\RequirePackage{times}
\RequirePackage{caption}
\RequirePackage{subcaption}
\RequirePackage{float}
\RequirePackage{algorithmicx}
\RequirePackage{algorithm}
\RequirePackage[noend]{algpseudocode}
\usepackage{array}
\usepackage{tabulary}
\newcolumntype{K}[1]{>{\centering\arraybackslash}p{#1}}
\RequirePackage{pgfplots}
\def\BState{\State\hskip-\ALG@thistlm}

\newcommand{\startcompact}[1]{\par\vspace{-0.75em}\begin{#1}%
\allowdisplaybreaks\ignorespaces}
\newcommand{\stopcompact}[1]{\end{#1}\ignorespaces}

\newenvironment{rcases}
  {\left.\begin{aligned}}
  {\end{aligned}\right\rbrace}
\RequirePackage[%
    font={small,sf},
    labelfont=bf,
    format=hang,    
    format=plain,
    margin=0pt,
    width=0.8\textwidth,
]{caption}
\tikzset{
    state/.style={
           rectangle,
           rounded corners,
           draw=black, very thick,
           minimum height=2em,
           inner sep=2pt,
           text centered,
           },
}
\usetikzlibrary{arrows,automata}
\usetikzlibrary{positioning}
\usetikzlibrary{arrows,shapes}
\RequirePackage{multirow}
\RequirePackage{hhline}
\RequirePackage{afterpage,lipsum}
\RequirePackage[labelformat=simple]{subcaption}

\newenvironment{customlegend}[1][]{%
\begingroup
\csname pgfplots@init@cleared@structures\endcsname
\pgfplotsset{#1}
}{%
\csname pgfplots@createlegend\endcsname
\endgroup
}%

\def\addlegendimage{\csname pgfplots@addlegendimage\endcsname}

\begin{document}


\title{Fast Solution of the Linearized Poisson-Boltzmann Equation with nonaffine Parametrized 
Boundary Conditions Using the Reduced Basis Method} 



\author{Cleophas Kweyu}
\email[]{kweyu@mpi-magdeburg.mpg.de}
\altaffiliation{Computational Methods in Systems and Control Theory (CSC)}
\affiliation{Max Planck Institute for Dynamics of Complex Technical Systems, Sandtorstr. 1, 39106, 
Magdeburg, Germany}

\author{Lihong Feng}
\email[]{feng@mpi-magdeburg.mpg.de}
\altaffiliation{Computational Methods in Systems and Control Theory (CSC)}
\affiliation{Max Planck Institute for Dynamics of Complex Technical Systems, Sandtorstr. 1, 39106, 
Magdeburg, Germany}

\author{Matthias Stein}
\email[]{matthias.stein@mpi-magdeburg.mpg.de}
\altaffiliation{Molecular Simulations and Design (MSD)}
\affiliation{Max Planck Institute for Dynamics of Complex Technical Systems, Sandtorstr. 1, 39106, 
Magdeburg, Germany}

\author{Peter Benner}
\email[]{benner@mpi-magdeburg.mpg.de}
\altaffiliation{Computational Methods in Systems and Control Theory (CSC)}
\affiliation{Max Planck Institute for Dynamics of Complex Technical Systems, Sandtorstr. 1, 39106, 
Magdeburg, Germany}

\date{\today}

\begin{abstract}
The Poisson-Boltzmann equation (PBE) is a nonlinear elliptic parametrized partial differential 
equation that arises in biomolecular modeling and is a fundamental tool for structural biology. It is 
used to calculate electrostatic potentials around an ensemble of fixed charges immersed in an ionic 
solution. It can also be used to estimate the electrostatic contribution to the free energy of a system. 
Efficient numerical computation of the PBE yields a high number of degrees of freedom in the resultant 
algebraic system of equations, ranging from several hundred thousands to millions. Coupled with the fact 
that in most cases the PBE requires to be solved multiple times for a large number of system configurations, 
this poses great computational challenges to conventional numerical techniques. To accelerate such 
computations, we here present the reduced basis method (RBM) which greatly reduces this computational 
complexity by constructing a reduced order model of typically low dimension. In this study, we employ a 
simple version of the PBE for proof of concept and discretize the linearized PBE (LPBE) with a centered 
finite difference scheme. The resultant linear system is solved by the aggregation-based algebraic multigrid  
(AGMG) method at different samples of ionic strength on a three-dimensional Cartesian grid. The discretized 
LPBE, which we call the high-fidelity full order model (FOM), yields solution as accurate as other LPBE 
solvers. We then apply the RBM to FOM. The discrete empirical interpolation method (DEIM) is applied 
to the Dirichlet boundary conditions which are nonaffine with the parameter (ionic strength), to reduce the 
complexity of the reduced order model (ROM). From the numerical results, we notice that the RBM reduces the 
model order from $\mathcal{N} = 2\times 10^{6}$ to $N = 6$ at an accuracy of $10^{-9}$ and reduces 
computational time by a factor of approximately $7,600$. DEIM, on the other hand, is also used in the 
offline-online phase of solving the ROM for different values of parameters which provides a speed-up of 
$20$ for a single iteration of the greedy algorithm.
\end{abstract}

\pacs{}

\keywords{Reduced basis method, Poisson-Boltzmann equation, ionic strength,
finite differences scheme, aggregation-based algebraic multigrid method, discrete empirical 
interpolation method.}

\maketitle 

\section{Introduction}
\subsection{Electrostatic Interactions in Biomolecular Systems}
Electrostatic interactions are important in biological processes such as molecular recognition, 
enzyme catalysis, and biomolecular encounter rates. A significant challenge in computational 
biology has been to model these interactions accurately and efficiently. This is because 
biomolecules are surrounded by solvent molecules and therefore, the solvent effects must be 
considered during modeling. There are two main groups of computational approaches which are used to 
model electrostatic interactions based on how the solvent is treated. Explicit methods place the 
solvent molecules around the biomolecule while implicit methods consider the solvent molecules as 
a continuum \cite{Wang2010,Fogolari2002}. 

The Poisson-Boltzmann equation (PBE) is one of the most popular implicit solvent models which 
describes the solvent in a continuum model through the Boltzmann distribution. The PBE solves 
the electrostatic potential in the entire domain which comprises both the molecule and the solvent. 
From this potential, further information can be obtained at various regions of interest and for 
different applications. Firstly, the electrostatic potential at the biomolecular surface, commonly 
known as electrostatic surface potential, can provide insights into possible docking sites for other 
small or large molecules. Secondly, the potential outside the biomolecule can provide information 
about the free energy of interaction of small molecules at different positions in the vicinity of 
the biomolecule. Thirdly, free energy of a biomolecule can be determined, which provides 
information about the molecule's stability. Finally, the electrostatic field can be estimated from 
which the mean atomic forces can be derived. More information can be found in 
\cite{Fogolari99,Fogolari2002,Honig95,SharpHonig90}.

Analytical solutions of the PBE are only possible under the assumption that the biomolecules of 
interest have regular shapes, for example, spheres or cylinders. And even if these solutions exist, 
they are still quite complex. However, these are not realistic because biomolecules have 
irregular shapes or geometries and charge distributions \cite{Holst94, Dong2008}. This makes it 
necessary to apply numerical techniques to the PBE and the first of such methods were introduced 
in \cite{Warwicker1982} where the electrostatic potential was determined at the active site of a 
protein (or enzyme). The most popular numerical techniques in this regard are based on 
discretization of the domain of interest into small regions and employ the finite difference 
methods (FDM) \cite{Baker2001, Wang2010}, the finite element methods (FEM) \cite{Baker2001, 
Holst2000}, or the boundary element methods (BEM) \cite{Boschitsch2004, Zhou1993}. A thorough 
review of the numerical methods for solving the PBE can be found in \cite{Lu2008}.  

All of the aforementioned numerical methods have one major advantage in common. It is possible to 
employ ``electrostatic focussing'', which enables users to apply relatively coarse grids for 
the entire calculations and very fine grids in regions of great interest such as the binding 
or active sites of macrobiomolecules. This adaptivity provides highly accurate local solutions 
to the PBE at reduced computational costs \cite{Baker2005}. However, the BEM has the drawback 
of being applicable only to the LPBE and thus limiting its general use. Numerous software 
packages have been developed to solve the PBE and some of the major ones include the adaptive 
Poisson-Boltzmann solver (APBS) \cite{Baker2001} and Delphi \cite{Rocchia2001}. There are also 
recent developments regarding the PBE theory which include, the treatment of the biomolecular 
system as an interface problem, the extensive studies on the nonlinear PBE, among others, see 
\Cref{PBE_theory} for more details.
 
Due to the limited computational memory and speed, solving the PBE efficiently is still 
computationally challenging and affecting the accuracy of the numerical solutions. This is 
due to the following reasons. Firstly, electrostatic interactions are long-ranged and therefore, 
the electrostatic potential decays exponentially over large distances, see equation
(\ref{eq:Decay_expo}). This requires an infinite domain which is infeasible in practice. Secondly, 
biomolecules of interest comprise of thousands to millions of atoms which require 
a large domain to accommodate both the biomolecule and the solvent. To circumvent these 
challenges, it is customary to choose a truncated domain of at least three times the size of the 
biomolecule so as to accurately approximate boundary conditions \cite{Holst94}. Nonetheless, this 
still leads to a very large algebraic system consisting of several hundreds of thousands to 
millions of degrees of freedom. It becomes even more difficult if the PBE is incorporated in a 
typical dynamics simulation which involves millions of time steps or in a multi-query task where 
the solution is solved many times for varying parameter values such as the ionic strength 
\cite{Wang2010}.

The computational complexity arising from the resultant high-dimensional system can be greatly 
reduced by applying model order reduction (MOR) techniques. The main goal of MOR is to construct 
a reduced-order model (ROM) of typically low dimension, whose solution retains all the important 
information of the high-fidelity system at a greatly reduced computational effort. Because the PBE 
is a parametrized PDE (PPDE), we apply the reduced basis method (RBM) which falls into the class of 
parametrized MOR (PMOR) techniques \cite{morBenGW15}. However, it is important to note that the RBM 
is not an independent numerical technique; hence its accuracy depends on that of the underlying technique 
which is used to discretize the PBE \cite{morBenGW15,morEft11}. In this paper, we discretize the PBE 
using FDM before applying the RBM.

The benefits of the RBM, or the MOR in general, become obvious when the same problem has to be 
solved for a large number of parameter values. In our study, the break-even point is about $10$, 
and thus, the RBM becomes very effective if dozens or more parameter configurations need to be 
evaluated.

\begin{figure}[t]
  \centering
    \includegraphics[width=5.5cm]{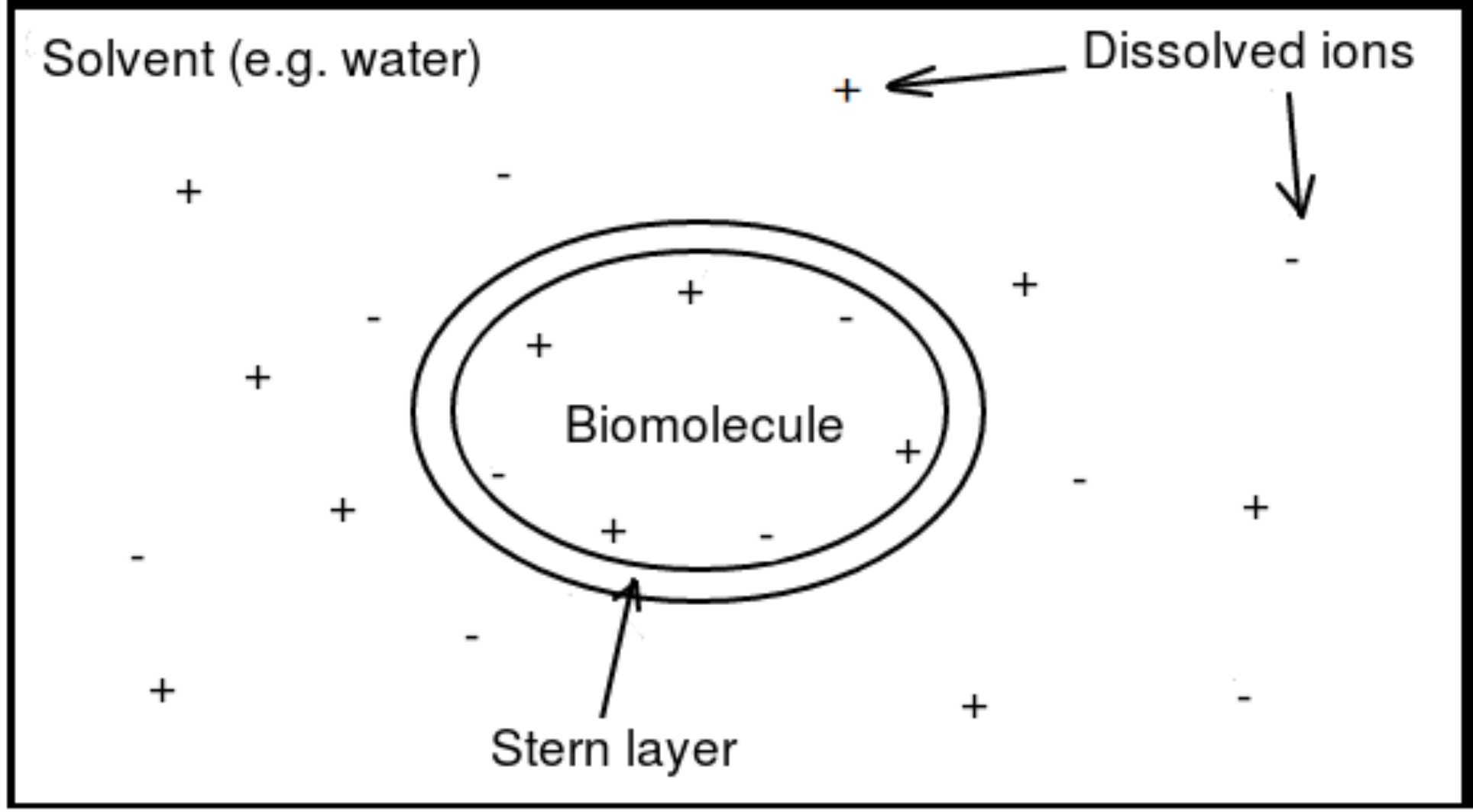}
    \caption{\label{fig:DB_model}
    2-D view of Debye-H\"uckel model.}
\end{figure}

Here, we consider a protein molecule immersed in ionic solution at physiological concentration, and 
determine the electrostatic potential triggered by the interaction between the two particles, see 
Figure \ref{fig:DB_model}. The electrolyte here is of monovalent type, implying that the ionic 
strength is equivalent to the concentration of the ions. The ionic strength is a physical parameter 
of the PBE, and we determine the electrostatic potential under variation of this parameter.
This paper is an extension of the ECCOMAS Congress 2016 proceedings paper, 
\cite{Kweyu2016} with the following additional key inputs. Firstly, we employ nonaffine Dirichlet 
boundary conditions given in \Cref{DH_conditions} to replace the zero Dirichlet boundary 
conditions in the former. Secondly, and as a consequence of the nonaffine parameter dependence 
of these boundary conditions, we apply the discrete empirical interpolation method (DEIM) to reduce 
the resultant complexity in the reduced order model (ROM) during the online phase of the reduced 
basis method (RBM), see \Cref{Complexity} and \Cref{DEIM} \cite{morChaS10}. Lastly, we apply finite 
volume discretization to the dielectric coefficient function instead of taking the averages of 
the dielectric values between two neighbouring grid points. This is meant to reduce the truncation 
error as explained in \Cref{Dielectric_kappa}.

\subsection{An overview of Poisson-Boltzmann Theory}\label{PBE_theory}

There are numerous ways and reviews on the derivation of the PBE. The simplest stems from the 
Poisson equation \cite{Jackson1975, Neves-Petersen2003} (in SI units),
\begin{equation}\label{eq:PE}
 -\vec{\nabla}.(\epsilon(x)\vec{\nabla}u(x)) = \rho(x), \quad \textrm{in} \quad \Omega \in \mathbb{R}^3,
\end{equation}
which describes the electrostatic potential $u(x)$ at a point $x \in \Omega$. The term $\rho(x)$ is 
the charge distribution which generates the potential in a region with a position-dependent and 
piecewise constant dielectric function $\epsilon(x)$. Equation (\ref{eq:PE}) is generally solved in a 
finite domain $\Omega$ subject to Dirichlet boundary conditions $u(x) = g(x)$ on 
$\partial{\Omega}$. Usually, $g(x)$ employs an analytic and asymptotically correct form of the 
electrostatic potential and therefore, the domain must be sufficiently large to ensure an accurate 
approximation of the boundary conditions \cite{Dong2008}.

To obtain the PBE from equation (\ref{eq:PE}), we consider two contributions to the charge 
distribution $\rho(x)$: the ``fixed'' solute charges $\rho_f(x)$ and the aqueous ``mobile'' ions in 
the solvent $\rho_m(x)$. The $N_m$ partial atomic point charges ($z_i$) of the biomolecule are 
modeled as a sum of delta distributions at each atomic center $x_i$, for $i = 1, \ldots, N_m$, that is,
\begin{equation}\label{eq:fixed_charges}
 \rho_f(x) = \frac{4\pi e^2}{k_B T}\sum_{i=1}^{N_m}z_i\delta(x-x_i).
\end{equation}
The term $e/{k_BT}$ is the scaling coefficient which ensures that the electrostatic potential 
is dimensionless, where $e$ is the electron charge and $k_BT$ is the thermal energy of the system and 
is comprised of the Boltzmann constant $k_B$ and the absolute temperature $T$. The total charge of each 
atom is $ez_i$. 

On the other hand, the solvent is modeled as a continuum through the Boltzmann distribution which 
leads to the mobile ion charge distribution
\begin{equation}\label{eq:mobile_ions}
 \rho_m(x) = \frac{4\pi e^2}{k_B T}\sum_{j=1}^{m}c_jq_je^{-q_ju(x)-V_j(x)},
\end{equation}
where we have $m$ mobile ion species with charges $q_j$ and bulk concentrations $c_j$. The term 
$V_j(x)$ is the steric potential which prevents an overlap between the biomolecule and the 
counterions. For monovalent electrolytes, whose ions are in a $1:1$ ratio, for example, 
$NaCl$, equation (\ref{eq:mobile_ions}) reduces to 
\begin{equation}\label{eq:mobile_ions2}
 \rho_m(x) = -\kappa^2(x)\sinh(u(x)),
\end{equation}
where the kappa function $\kappa^2(x)$ is position-dependent and piecewise constant; it describes both 
the ion accessibility through $e^{-V(x)}$ and the bulk ionic strength (or concentration) 
\cite{Baker2005}.

We eventually obtain the PBE by combining the two expressions for the charge distributions in 
(\ref{eq:fixed_charges}) and (\ref{eq:mobile_ions2}) with the Poisson equation (\ref{eq:PE}) 
for a monovalent electrolyte, 
\begin{equation}\label{eq:PBE}
 -\vec{\nabla}.(\epsilon(x)\vec{\nabla}u(x)) + \bar{\kappa}^2(x) \sinh(u(x)) = 
 \frac{4\pi e^2}{k_B T}\sum_{i=1}^{N_m}z_i\delta(x-x_i), \quad \textrm{in} \quad \Omega \in 
 \mathbb{R}^3,
\end{equation}
subject to
\begin{equation}\label{eq:Dirichlet}
 u(x) = g(x) \quad  \textrm{on} \quad \partial{\Omega},     
\end{equation}
where
\begin{equation}\label{eq:Decay_expo}
u(\infty) = 0 \, \implies u(x_{\max})\rightarrow 0 \,\,\textrm{as} \,\, 
\left|x_{\max}\right|\rightarrow \infty.    
\end{equation}
In equation (\ref{eq:PBE}), $ u(x) = {e\psi(x)}/{k_B T}$ is the dimensionless potential 
scaled by $e/{k_B T}$ and $\psi(x)$ is the original electrostatic potential in 
centimeter-gram-second (cgs) units at $ x \in \mathbb{R}^3 $. The terms $\epsilon(x)$ and 
$\bar{k}^2(x)$ are discontinuous functions at the interface between the charged biomolecule and 
the solvent, and at an ion exclusion region (Stern layer) surrounding the molecule, respectively. 
The term $k^2 = {8\pi e^2 I}/{1000\epsilon k_BT}$ is a function of the ionic strength 
$I = 1/2\sum_{i=1}^{N}c_iz_i^2$. The function $g(x)$ represents the Dirichlet 
boundary conditions which are discussed in detail in \Cref{DH_conditions} and are nonaffine in 
the parameter $I$. Equation (\ref{eq:Decay_expo}) shows that the electrostatic potential decays to 
zero exponentially as the position approaches infinity. Details on mapping $\epsilon(x)$ and 
$\bar{k}^2(x)$ onto a computational grid can be found in \cite{Baker2001}. The PBE (\ref{eq:PBE}) 
poses severe computational challenges in both analytical and numerical approaches due to the 
infinite (unbounded) domain in (\ref{eq:Decay_expo}), delta distributions, rapid nonlinearity, and 
discontinuous coefficients \cite{Holst94, Holst2000}.  

The PBE (\ref{eq:PBE}) can be linearized under the assumption that the electrostatic potential 
is very small relative to the thermal energy $k_BT$ \cite{Fogolari2002}. 
Therefore, the nonlinear function $\textrm{sinh}(u(x))$ can be expanded into a Taylor series 
\begin{equation}\label{eq:Taylor}
  \sinh(u(x)) = u(x) + \frac{u(x)^3}{3!} + \frac{u(x)^5}{5!} + \ldots,
\end{equation}
and only the first term is retained. We obtain the linearized PBE (LPBE) given by
\begin{equation}\label{eq:LPBE}
   -\vec{\nabla}.(\epsilon(x)\vec{\nabla}u(x)) + \bar{k}^2(x)u(x) = 
   (\frac{4\pi e^2}{k_B T})\sum_{i=1}^{N_m}z_i\delta(x-x_i).
\end{equation}
Usually, proteins are not highly charged, and it suffices to consider the linearized PBE (LPBE). One can 
still obtain accurate results because the higher order terms in (\ref{eq:Taylor}) do not provide 
a significant contribution. However, we must note that the LPBE can give inaccurate 
results for highly charged biomolecules such as the DNA and RNA (nucleic acids), phospholipid 
membranes, and polylysine \cite{Honig95}. More information about the PBE, including its 
derivation from first principles, can be found in \cite{Holst94}.

It is worth noting that there are recent developments of the PBE theory. Firstly, the biomolecular 
system has been considered as an interface problem which requires solution decomposition techniques to 
get rid of the solution singularities caused by the Dirac-delta distributions on the right hand 
side of (\ref{eq:LPBE}) or (\ref{eq:PBE}). This has been discussed, for example, in  
\cite{Chen2009, Xie2014, Xie2016} where the LPBE has been modified into the form
\begin{equation}\label{eq:LPBE_interface}
\begin{rcases}
\begin{aligned}
 -\epsilon_p \Delta u(x) = \alpha \sum_{i=1}^{N_m}z_i\delta(x-x_i),  \quad & x \in D_p,\\
 -\epsilon_s \Delta u(x) + \kappa^2 u(x) = 0, \quad & x \in D_s,\\
 u(s^+) = u(s^-), \quad \epsilon_s \frac{\partial u(s^+)}{\partial n(s)} 
 = \epsilon_p \frac{\partial u(s^-)}{\partial n(s)}, \quad & s \in \varGamma,\\
 u(s) = g(s), \quad & s \in \partial \Omega, \quad
\end{aligned}
\end{rcases}
\end{equation}
where $\alpha$ is a constant, $D_p$ the protein domain, $D_s$ the solvent domain and $\varGamma$ the 
interface between the protein and the solvent. The PBE (nonlinear) has also been extensively solved 
as an interface problem \cite{Chen2009, Xie2014}.

The interface problem in (\ref{eq:LPBE_interface}) is more accurate than the model in (\ref{eq:LPBE}) 
considered in this study, because the local or short-range potentials generated by the 
Dirac-delta distributions are computed independent of the long-range potentials, thus avoiding 
errors. However, this model is still computationally expensive because the numerical calculations by 
conventional methods are in $\mathcal{O}(\mathcal{N}^3)$, (commonly known as the ``curse of 
dimensionality''), where $\mathcal{N}$ is the dimension of the system in one direction. Therefore, we use 
the simple model (\ref{eq:LPBE}) for the purpose of introducing and validating the RBM. Considering the 
interface problem would be our next step. 

Secondly, studies on a variational problem of minimizing a mean-field variational electrostatic free-energy 
functional have been conducted \cite{LiWeSh:16}. This has been done in order to investigate the dependence of 
dielectric coefficient on local ionic concentrations and its effect on the equilibrium properties of 
electrostatic interactions in an ionic solution which was proposed, for instance, in \cite{SaGuMa:10}. Results 
show that indeed the dielectric coefficient depends on the local ionic concentrations and this dependence can 
be expressed as a mathematical function which is continuous, monotonically decreasing, and convex 
\cite{LiWeSh:16}.

\subsection{Applications and Post-processing of the PBE solution}

The resultant electrostatic potential for the entire system can be used to calculate 
electrostatic free energies and electrostatic forces. The electrostatic free energy represents the 
work needed to assemble the biomolecule and is obtained by integration of the potential over a 
given domain of interest \cite{Vergara-Perez2016, Dong2008}. For the LPBE, this energy is given by
\begin{equation}\label{eq:Elec_energy}
G_{elec}\lbrack u(x)\rbrack = \frac{1}{2}\int_{\Omega}\rho_f u(x)dx = 
\frac{1}{2}\sum_{i=1}^{N}z_iu(x_i),
\end{equation}
where $u(x_i)$ is the mean electrostatic potential acting on an atom $i$ located at position 
$x_i$ and carrying a charge $z_i$. The integral in (\ref{eq:Elec_energy}) can be seen as the 
integral of polarization energy which is equivalent to the sum of interactions between 
charges and their respective potentials. On the other hand, it is also possible to differentiate 
the energy functional in (\ref{eq:Elec_energy}) with respect to atomic positions to obtain the 
electrostatic force on each atom \cite{Gilson1993, Baker2005, Dong2008}.

\begin{figure}[H]
  \centering
  \begin{minipage}[b]{0.19\textwidth}
    \includegraphics[width=\textwidth]{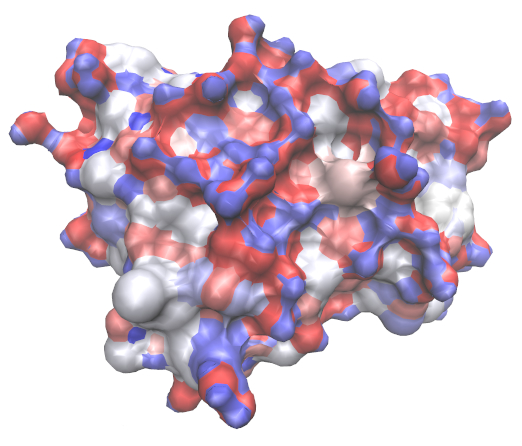}
  \end{minipage}
  \quad \quad
  \begin{minipage}[b]{0.19\textwidth}
    \includegraphics[width=\textwidth]{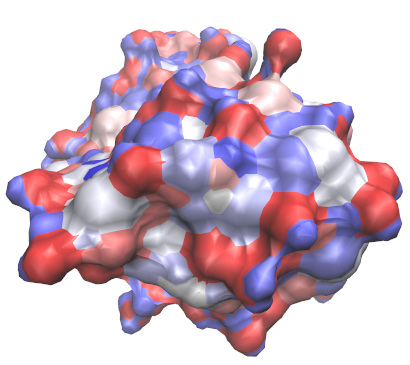}
  \end{minipage}
  \caption{\label{fig:Surf_pot}%
 Electrostatic potential mapped onto protein \\surface of \textit{fasciculin 1} toxin CPDB entry 1FAS.}
\end{figure}

The electrostatic potential can also be evaluated on the surface of the biomolecule (electrostatic 
surface potential). It is used to provide information about the interaction between the biomolecule and 
other biomolecules or ligands or ions in its vicinity. Figure \ref{fig:Surf_pot} shows the electrostatic 
potential mapped onto the surface of the protein \textit{fasciculin 1} and was generated by the Visual 
Molecular Dynamics (VMD) software at different orientations \cite{Humphrey1996}. The electrostatic potential 
is computed by our FDM solver. The red colour represents regions of negative potential, the blue colour 
represents regions of positive potential, and the white colour represents neutral regions. 

The outline of this paper is as follows: In \Cref{FDM_discretn}, we provide a glimpse on the finite difference 
discretization of the LPBE and those of the dielectric coefficient and kappa functions, charge 
densities, as well as their respective mappings to the computational grid. In \Cref{RBM_summary}, we 
provide the basics of the RBM which include the problem formulation, the solution manifold, the 
greedy algorithm, the discrete empirical interpolation method (DEIM), and the 
\textit{a posteriori} error estimation. In \Cref{Numerical_results}, we provide numerical results of the FOM 
(via the FDM) and those of the ROM (via RBM and DEIM). Conclusions and some ideas on future work are 
given in the end. 

\section{Discretization of the Poisson-Boltzmann Equation}\label{FDM_discretn}

\subsection{Finite Difference Discretization}

We discretize the LPBE in (2) with a centered finite differences scheme to obtain the algebraic 
linear system as below,
\startcompact{small}
\begin{multline}\label{eq:FD_eqn}
 \displaystyle -\frac{H}{dx^2}\epsilon_{i+\frac{1}{2},j,k}^x(u_{i+1,j,k}-u_{i,j,k})+
 \frac{H}{dx^2}\epsilon_{i-\frac{1}{2},j,k}^x(u_{i,j,k}-u_{i-1,j,k})\\
 -\frac{H}{dy^2}\epsilon_{i,j+\frac{1}{2},k}^y(u_{i,j+1,k}-u_{i,j,k})+ 
 \frac{H}{dy^2}\epsilon_{i,j-\frac{1}{2},k}^y(u_{i,j,k}-u_{i,j-1,k})\\
 -\frac{H}{dz^2}\epsilon_{i,j,k+\frac{1}{2}}^z(u_{i,j,k+1}-u_{i,j,k})+
 \frac{H}{dz^2}\epsilon_{i,j,k-\frac{1}{2}}^z(u_{i,j,k}-u_{i,j,k-1})\\
 + H\bar{\kappa}^2_{i,j,k}u_{i,j,k} = HCq_{i,j,k},
\end{multline}
\stopcompact{small}
where $H = dx\times dy\times dz$ is a scaling factor, $q_{i,j,k}$ is the discretized molecular 
charge density and $C = {4\pi e^2}/{k_B T}$. 
It is important to choose efficient algorithms and parameters to be used in the discretization of 
the charge density distribution, the kappa, and the dielectric functions that appear in the LPBE 
for the accuracy of the mean electrostatic potential solution. An efficient method is usually 
chosen to partition the domain into regions of solute (or biomolecule) and the solvent 
dielectric. Some of the key methods employed in APBS are the molecular surface and cubic-spline 
surface methods \cite{Vergara-Perez2016}. In the following subsections, we provide some insights 
into these discretizations.

\subsection{Calculation of Dielectric Constant Distribution and Kappa function}\label{Dielectric_kappa}

We notice that the dielectric constant $\epsilon$ in equation (\ref{eq:FD_eqn}), is 
discretized at half grid, and therefore, we use a staggered mesh which results in three arrays 
(in $x$, $y$, and $z$ directions) representing the shifted dielectric values on different grids. 
This intends to fully take advantage of the finite volume discretization in order to minimize the 
solution error by increasing the spatial resolution. The dielectric coefficients and 
kappa functions which are piecewise constant, are mapped according to the following conditions,
\begin{eqnarray}
 \epsilon(x) =
  \begin{cases}
   2 & \text{if } x \in \Omega_1\\
   78.54 & \text{if } x \in \Omega_2 \,\, \text{or} \,\, \Omega_3
  \end{cases}, \\
 \bar{k}(x) =
  \begin{cases}
   0 & \text{if } x \in \Omega_1 \,\, \text{or} \,\, \Omega_2\\
    \sqrt{\epsilon_3}k & \text{if } x \in \Omega_3 
  \end{cases},
\end{eqnarray}
where $\Omega_1$ is the region occupied by the protein molecule, $\Omega_2$ is the ion-exclusion 
layer, and $\Omega_3$ is the region occupied by the ionic solution.

Techniques used to map the dielectric and kappa functions onto the grid include, among others, the 
molecular surface, and the smoothed molecular surface, which are calculated using the Connolly approach 
\cite{Connolly83} and the cubic-spline surface. For more information see \cite{Vergara-Perez2016}. 
The cubic-spline surface method, which is our method of choice, is more suitable than the other two 
because it is possible to evaluate the gradient of the mean electrostatic potential such as in the 
determination of the solvated or polar forces. This method introduces an intermediate dielectric region 
at the interface between the solute and the solvent because the kappa and dielectric maps are built on 
a cubic-spline surface. This smoothes the transition of the functions to circumvent discontinuities 
inherent in them \cite{Baker2001, Vergara-Perez2016}.

\subsection{Calculation of Charge Densities}

The molecular charge density (right-hand side of the LPBE (\ref{eq:LPBE})) can be obtained from any 
file with atomic coordinates, charges, and radii. However, these atomic coordinates may not 
coincide with any of our grid points. Therefore, it is necessary to find an efficient method of 
spreading the point charges (summation term in LPBE) to the grid points. Several methods are 
available to map or spread the charges onto the grid points, e.g. in the 
APBS software package. Trilinear interpolation (or linear spline) in which charges are mapped 
onto nearest-neighbour grids, results in potentials which are very sensitive to the grid 
resolution. Cubic B-spline interpolation where charges are mapped to two layers of grid points, 
has an average sensitivity to the grid setup, and quintic B-spline interpolation has the lowest 
sensitivity to grid spacing because charges are spread out to three layers of the grid points 
\cite{Baker2001}. 

In this study, we use the cubic B-spline interpolation (basis spline) method which maps the charges 
to the nearest and next-nearest grid points. Although computationally expensive, this method provides 
softer or smoother distributions of charges which subsequently reduces the sensitivity of the mean 
electrostatic potential solutions to the grid spacing \cite{Vergara-Perez2016}.

\subsection{Dirichlet Boundary Conditions}\label{DH_conditions}

Analytical solutions to the LPBE can only be obtained for systems with simple geometries, for 
example, spherical and cylindrical systems. Equation (\ref{eq:DH_soln}) shows an analytical 
solution for a spherical molecule with uniform charge (Born ion) \cite{Holst94}. From this equation, 
we can obtain two different kinds of Dirichlet boundary conditions, the single Debye-H\"uckel (SDH) 
and multiple Debye-H\"uckel (MDH). For the former, we assume that all the atomic charges are 
collected into a single charge at the center of the solute approximated by a sphere. This kind of 
boundary condition is suitable when the boundary is sufficiently far from the biomolecule. On the 
other hand, the latter assumes the superposition of the contribution of each atomic charge (i.e. 
multiple, non-interacting spheres with point charges) with respective radius. This kind of boundary 
condition is more accurate than SDH for closer boundaries but can be computationally expensive for 
large biomolecules. 

In this study, we employ the MDH type \cite{Baker2001, Rocchia2005},
\begin{equation}\label{eq:DH_soln}
 u(x) = (\frac{e^2}{k_B T}) \sum\limits_{i=1}^{Nm}\frac{z_ie^{-\kappa(d-a_i)}}{\epsilon_w (1+ka_i)d} \quad
 \textrm{on} \quad \partial{\Omega}, \quad d = \lvert x-x_i \rvert.      
\end{equation}
Here, $z_i$ are the point partial charges of the protein, $\epsilon_w$ is the solvent dielectric, 
$\kappa = \bar{\kappa}/ \sqrt{\epsilon_w}$ is a function of the ionic strength of the solution, $a_i$ are 
the atomic radii, and $N_m$ is the total number of point partial charges in the protein.

\section{Essentials of the Reduced Basis Method}\label{RBM_summary}

The Reduced basis method (RBM) and proper orthogonal decomposition (POD) are examples of popular 
projection-based parametrized model order reduction (PMOR) techniques. The main goal of these 
techniques is to generate a parametric ROM which accurately approximates 
the original full order model (FOM) of high dimension over varying parameter values 
\cite{morBenGW15, morEft11, morRozHP08}. The RBM exploits an offline/online procedure which ensures an 
accurate approximation of the high-fidelity solution in a rapid and inexpensive manner and is widely 
applicable in real-time and many-query scenarios. For a thorough review, see \cite{morBenGW15}. 

We consider a physical domain $\Omega \subset \mathbb{R}^3$ with boundary 
$\partial \Omega$, and a parameter domain $D \subset \mathbb{R}$. The LPBE (\ref{eq:LPBE}) is 
discretized with the centered finite difference scheme (\ref{eq:FD_eqn}) on $\Omega$ and Dirichlet 
boundary conditions (\ref{eq:Dirichlet}) obtained from (\ref{eq:DH_soln}) are applied. The resultant 
discrete problem of the LPBE becomes, for any $\mu \in D$, find $u^{\mathcal{N}}(\mu)$ that 
satisfies the linear system
\begin{equation}\label{eq:FOM}
 A(\mu)u^{\mathcal{N}}(\mu) = f(\mu),
 \quad \mu \in D,
\end{equation}
where $A(\mu) \in \mathbb{R}^{\mathcal{N}\times \mathcal{N}}$ 
and $f(\mu) \in \mathbb{R}^{\mathcal{N}}$. The matrix $A(\mu)$ can also be written as a 
parameter-affine matrix, 
\begin{equation}
 A(\mu) = \sum\limits_{i=1}^{Q}\Theta_i(\mu)A_i,
\end{equation}
where $Q\in \mathbb{N}$, $\Theta_i$ are scalar coefficient functions, and $A_i$ are the parameter 
independent matrices.
The $\mathcal{N}\times \mathcal{N}$ system is indeed computationally expensive to be solved for an 
accurate approximation of $u(\mu)$ because the dimension $\mathcal{N}$ is approximately 
$2\times10^6$ in our problem. Therefore, we apply the RBM to save computational costs by providing an 
accurate approximation of $u^{\mathcal{N}}(\mu)$ at a greatly reduced dimension of
$N \ll \mathcal{N}$. The ROM is given by (\ref{eq:ROM}).

However, as detailed in \Cref{Complexity}, we encounter some computational complexity in the 
online phase of RBM which is caused by the nonaffine parameter dependence in the right-hand side 
vector $f(\mu)$ from the boundary condition (\ref{eq:DH_soln}). The parameter, the ionic strength, 
resides in the kappa term $\kappa$ in the exponential function. This violates one of the key 
assumptions of the RBM which requires that all the system matrices and vectors must be 
affinely dependent on the parameter so that the offline/online decomposition is natural 
\cite{morHest16}. To circumvent this problem, we propose to apply an empirical interpolation method 
to reduce the complexity of the the online phase by avoiding the high-dimensional computation related to 
the vector $f(\mu)$. We provide some details in \Cref{DEIM}.

\subsection{The Solution Manifold and the Greedy Algorithm}

Another key assumption in RBM besides the affine parameter dependence, is the existence of a 
typically smooth and very low dimensional solution manifold which almost covers all the 
high-fidelity solutions of (\ref{eq:FOM}) under variation of parameters \cite{morEft11},
\begin{equation}\label{Manifold}
 \mathcal{M}^{\mathcal{N}} = \{ u^{\mathcal{N}}(\mu) : \mu \in D\}.
\end{equation}
The RB approximation space is then built upon this solution manifold and is given by the subspace 
spanned by the snapshots of the FOM. In other words, it is the subspace spanned by the high-fidelity 
$u^{\mathcal{N}}(\mu)$ solutions corresponding to a number of samples of the parameters, that is,
\begin{equation}\label{RB_space}
 \textrm{range}(V) = \textrm{span}\{ u^{\mathcal{N}}(\mu_1), ..., u^{\mathcal{N}}(\mu_l)\}, 
 \quad \forall \mu_1, ..., \mu_l \in D.
\end{equation}
The greedy algorithm as given in \Cref{euclid} is used to generate the reduced basis space 
(\ref{RB_space}) through an iterative procedure where a new basis is computed at each iteration 
\cite{morHesB13}. The RB space can be thought of being nested or hierarchical such that the 
previous basis set is a subset of the next and so on.
\begin{algorithm}[H]
  \caption{Greedy algorithm\cite{Hesthaven2014}}\label{euclid}
    \begin{algorithmic}[1]
    \Require A training set $\Xi \subset D$ including samples of $\mu$ covering the parameter domain 
	    $D$, i.e., $\Xi:= \{\mu_1, \ldots, \mu_l\}$. 
    \Ensure RB basis represented by the projection matrix $V$.
    \State Choose $\mu^* \in \Xi$ arbitrarily.
    \State Solve FOM (\ref{eq:FOM}) for $u^{\mathcal{N}}(\mu^*)$.
    \State $V_1 = [u^{\mathcal{N}}(\mu^*)]$, $N = 1$.
    \While{$\max\limits_{\mu\in \Xi}\Delta_N(\mu) \geq \epsilon$}
    \State $\mu^* = \textrm{arg}\max\limits_{\mu\in \Xi}\Delta_N(\mu)$.
    \State Solve FOM (\ref{eq:FOM}) for $u^{\mathcal{N}}(\mu^*)$.
    \State $V_{N+1} = [V_N \quad u^{\mathcal{N}}(\mu^*)]$.
    \State Orthonormalize the columns of $V_{N+1}$.
    \State $N = N + 1$.
    \EndWhile
    \State \textbf{end while}
  \end{algorithmic}
\end{algorithm}
The RB approximation is then formulated as, for any given $\mu \in D$, find $u_N(\mu) \in X_N$ 
which satisfies
\begin{equation}\label{eq:ROM}
 A_N(\mu)u_N(\mu) = f_N(\mu),
\end{equation}
where $A_N = V^TAV$ and $f_N(\mu) = V^Tf(\mu)$. $V$ is the orthonormal matrix computed from the greedy 
algorithm. From the fact that $N \ll \mathcal{N}$, solving the small dimensional reduced order 
model (ROM) is much cheaper than solving the high-fidelity model, the FOM (\ref{eq:FOM}) 
\cite{morEft11}. However, one problem still remains when computing the ROM. The computational 
complexity of evaluating the nonaffine function $f_N(\mu)$ still depends on the dimension of the FOM, 
as illustrated in \Cref{Complexity}. Efficient implementation of \Cref{euclid} depends on an efficient 
error estimation $\Delta_N(\mu)$ of the ROM, which is discussed in \Cref{Error_estmn}.

\subsection{Computational Complexity of the Reduced Order Model (ROM)}\label{Complexity}

To demystify the issue of computational complexity in the ROM, we can first rewrite (\ref{eq:FOM}) 
explicitly to illustrate the affine parameter decomposition on the left-hand side and the 
nonaffine right-hand side,
\begin{equation}\label{eq:Affine_A}
 (A_1 +\mu A_2)u^{\mathcal{N}} = \rho + b(\mu),
 \quad \mu \in D,
\end{equation}
where the matrix $A_1$ comes from the Laplacian operator term, $A_2$ is a diagonal matrix from the 
$kappa$ term, $\rho$ represents the charge density term and $b(\mu)$, the boundary conditions obtained 
from the analytical solution in (\ref{eq:DH_soln}). We can clearly notice the affine parameter 
decomposition of the matrix $A$ in (\ref{eq:FOM}) into $A_1$ and $\mu A_2$ in (\ref{eq:Affine_A}). 
However, the right-hand side function $b(\mu)$ is nonaffine in the parameter and therefore it cannot 
be decomposed in such a manner. Consider the ROM which is obtained by the greedy algorithm 
approach in \Cref{euclid} and a Galerkin projection,
\begin{equation}\label{eq:ROM_matrices}
(\underbrace{\hat{A_1}}_{N \times N} + 
 \mu \underbrace{\hat{A_2}}_{N \times N})\underbrace{u^N}_{N \times 1} = \underbrace{\hat{\rho}}_{N \times 1} 
 + \underbrace{V^T}_{N \times \mathcal{N}}\underbrace{b(\mu)}_{\mathcal{N} \times 1},
\end{equation}
where $\hat{A_1} = V^TA_1V$, $\hat{A_2} = V^TA_2V$, $\hat{\rho} = V^T\rho$, and 
$N \ll \mathcal{N}$.

It is clear from (\ref{eq:ROM_matrices}) that the last term of the right-hand side (RHS) 
still depends on the dimension $\mathcal{N}$ of the FOM while all the other matrices and vectors 
depend only on the dimension $N$ of the ROM, with $N \ll \mathcal{N}$. Therefore, the reduced order 
matrices on the left-hand side and the first vector on the right-hand side of 
(\ref{eq:ROM_matrices}) can be precomputed and stored during the offline phase, thereby providing 
a lot of computational savings. However, the term $V^Tb(\mu)$ cannot be precomputed because of 
the aforementioned nonaffine parameter dependence and therefore, the Galerkin projection involving 
matrix-vector products which are dependent on the dimension $\mathcal{N}$, has to be computed in 
the online phase of solving the ROM. 

In principle, we require $\mathcal{O}(2\mathcal{N}N)$ flops for these matrix-vector products and a full 
evaluation of the nonaffine analytical function (\ref{eq:DH_soln}) to obtain $V^Tb(\mu)$. This can be 
computationally expensive for a large $\mathcal{N}$, especially during the a posteriori error estimation 
(computing $\Delta_N (\mu)$), where the residual is computed $l$ times for varying parameter values 
$\mu_i$, $i = 1, \ldots, l$ for a single iteration of the greedy algorithm. The discrete empirical 
interpolation method (DEIM) is an approach to circumvent this problem in order to reduce the computational 
complexity of the nonaffine function. We discuss this technique at length in the next Subsection.

\subsection{Discrete Empirical Interpolation Method (DEIM)}\label{DEIM}

DEIM is a complexity reduction technique that was proposed in \cite{morChaS10} to overcome 
the drawback of the proper orthogonal decomposition (POD) approach for approximating a nonaffine 
(or nonlinear) parametrized function in the ROM during the online phase. This drawback is in the sense 
that the evaluation of the nonlinear/nonaffine function is still equivalent to that of computing the 
counterpart of the original system, which yields no computational savings. Therefore, the main idea of 
DEIM is to interpolate the nonlinear/nonaffine function by computing only a few entries of it, which 
dramatically reduces the computational complexity \cite{morChaS10, morWirSH14}. 

We provide a brief overview on using the singular value decomposition (SVD) to obtain the interpolation 
basis vectors. Firstly, we compute snapshots of the function $b(\mu)$ at a set of parameter $\mu$ in the 
training set $\Xi=\{\mu_1, \ldots, \mu_l\}\subset D$ and construct the snapshot matrix,
\begin{equation}
 F=[b(\mu_1), \ldots, b(\mu_l)] \in \mathbb{R}^{\mathcal{N}\times l}.
\end{equation}
Secondly, we compute its singular value decomposition (SVD),
\begin{equation}\label{eq:svd}
 F = U_F\Sigma W^T,
\end{equation}
where $U_F \in \mathbb{R}^{\mathcal{N}\times l}$, $\Sigma \in \mathbb{R}^{l\times l}$, and 
$W \in \mathbb{R}^{l\times l}$. Note that the matrices $U_F$ and $W$ are orthogonal, that is, 
$(U_F)^TU_F=W^TW = I_l$, $I_l \in \mathbb{R}^{l \times l}$ and $\Sigma = \textrm{diag}(\sigma_1, 
\ldots, \sigma_l)$, with $\sigma_1 \geq \ldots \geq \sigma_l \geq 0$. 

\begin{figure}[t]
  \centering
  \begin{subfigure}[b]{0.22\textwidth}
\includegraphics[width=\textwidth]{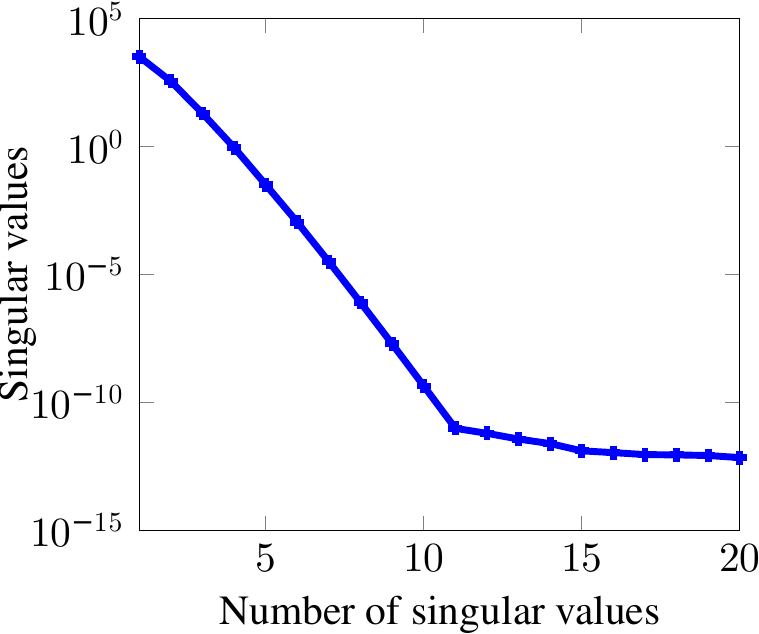}
\caption{20 singular values}
\label{fig:Sing_vals_20}
\end{subfigure}
\quad
\begin{subfigure}[b]{0.22\textwidth}
\includegraphics[width=\textwidth]{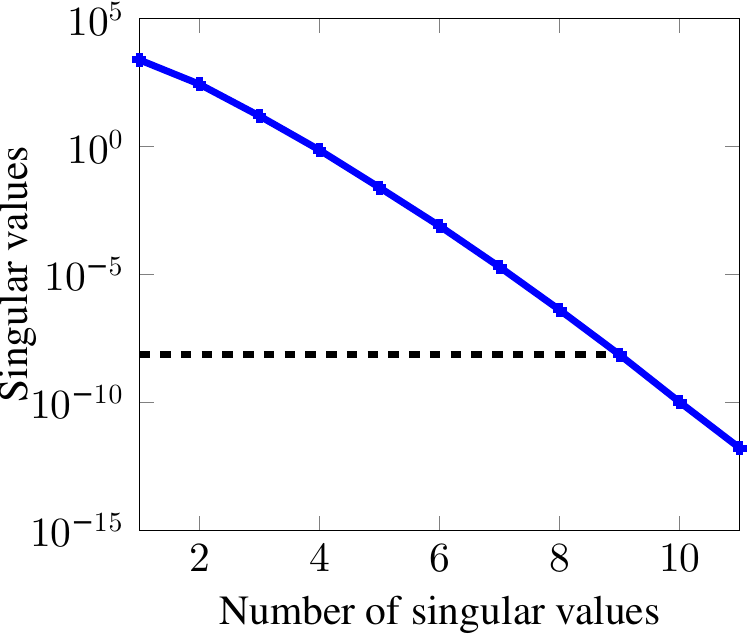}
\caption{11 singular values}
\label{fig:Sing_vals_11}
\end{subfigure}
  \caption{
  \label{fig:Sing_vals} %
  Decay of singular values of $\Sigma$ in (\ref{eq:svd}).}
\end{figure}

\Cref{fig:Sing_vals} shows the decay of the singular values of $\Sigma$ for the protein \textit{fasciculin 1}. 
\Cref{fig:Sing_vals_20} shows the behaviour of $20$ singular values with almost no decay from the $11th$ 
singular value. We discard these non-decaying singular values to obtain those in \Cref{fig:Sing_vals_11}. 
From the latter, we can actually truncate the singular values by selecting only the largest of them 
represented by $r \in \{ 1, \ldots, l\}$ that correspond to some required degree of accuracy. In this case, 
$l=11$ and $r=9$ which corresponds to an accuracy of $\epsilon_{svd}=O(10^{-10})$ in (\ref{eq:Error_svd}). 
The number $r$ plays an important role to select a basis set $\{u_i^F\}_{i=1}^r$ of rank $r$ from $U_F$ 
which solves the minimization problem \cite{morVol13},
\begin{equation}
\textrm{arg}\min\limits_{\{\tilde u_i\}_{i=1}^r}\sum_{j=1}^{l}\Arrowvert F_j -\sum_{i=1}^{r}\langle F_j,
\tilde u_i \rangle \tilde u_i\Arrowvert_2^2,  \quad \text{s.t.} \, \langle \tilde u_i,\tilde u_j \rangle 
= \delta_{ij},
\end{equation}
where $F_j$ is the $j$th column of the snapshot matrix $F$, and $\delta_{ij}$ is the usual Kronecker 
delta. 

The following criterion is used to truncate the largest singular values from \Cref{fig:Sing_vals} 
based on some desired accuracy, $\epsilon_{svd}$.
\startcompact{small}
\begin{equation}\label{eq:Error_svd}
 \frac{\sum\limits_{i=r+1}^{l}\sigma_i}{\sum\limits_{1=1}^{l}\sigma_i}<\epsilon_{svd},
\end{equation}
\stopcompact{small}
where $\sigma_i, i = 1, \ldots, l$ are the nonzero singular values of $F$. The dotted horizontal 
black line corresponds to $r = 9$ singular values and the corresponding singular vectors 
$\{u_i^F\}_{i=1}^r$ are used in the DEIM approximation.

DEIM overcomes the problem mentioned in \Cref{Complexity} by determining an interpolation of the 
nonaffine function $b(\mu)$. This is realized by approximating $b(\mu)$ with the linear combination of 
the basis vectors $U_F=[u_1^F, \ldots, u_r^F] \in \mathbb{R}^{\mathcal{N} \times r}$, i.e.
\begin{equation}
 b(\mu) \approx U_Fc(\mu),
\end{equation}
where $c(\mu)\in \mathbb{R}^r$ is the corresponding coefficient vector, and can be determined by 
assuming that $U_Fc(\mu)$ interpolates $b(\mu)$ at $r$ selected interpolation points, then,
\begin{equation}\label{eq:Overdet_sys}
 P^Tb(\mu) = P^TU_Fc(\mu),
\end{equation}
where $P$ is an index matrix given by
\begin{equation}
 P = [e_{\wp_1}, \ldots, e_{\wp_r}] \in \mathbb{R}^{\mathcal{N} \times r},
\end{equation} 
which consists of unit vectors $e_{\wp_i}$, $i = 1, \ldots, r$, where the indices $\wp_i$, are the DEIM 
interpolation points which are selected iteratively with a greedy algorithm. Suppose that $P^TU_F \in 
\mathbb{R}^{r \times r}$ is nondegenerate, then $c(\mu)$ can be determined from 
(\ref{eq:Overdet_sys}) by
\begin{equation}
c(\mu) = (P^TU_F)^{-1}P^Tb(\mu).
\end{equation}
Therefore, the function $b(\mu)$ in (14) can be approximated as
\begin{equation}\label{eq:DEIM_appro}
 b(\mu) \approx U_Fc(\mu) = U_F(P^TU_F)^{-1}P^Tb(\mu),
\end{equation}
so that the ROM in (\ref{eq:ROM_matrices}) with DEIM approximation becomes,
\begin{equation}
 (\underbrace{\hat{A_1}}_{N \times N} + 
\mu \underbrace{\hat{A_2}}_{N \times N})\underbrace{u^N(\mu)}_{N \times 1} = \underbrace{\hat{\rho}}_{N \times 1} 
 + \underbrace{V^TU_F(P^TU_F)^{-1}}_{N \times r}\underbrace{P^Tb(\mu)}_{r \times 1}. 
 \end{equation}
The interpolant $V^TU_F(P^TU_F)^{-1}P^Tb(\mu)$ can be computed a lot cheaper than $V^Tb(\mu)$ 
because we can precompute $V^TU_F(P^TU_F)^{-1}$ independent of the parameter $\mu$. Alternatively, we 
can also compute only those entries in $b(\mu)$ that correspond to the interpolation indices 
$\wp_i, i = 1, \ldots, r$, $r \ll \mathcal{N}$, i.e., $P^T b(\mu)$ instead of the entire $\mathcal{N}$ entries 
in $b(\mu)$. 

For the actual numerical implementation of the interpolation (\ref{eq:DEIM_appro}), the 
matrix $P$ needs not be explicitly applied. Instead, only the interpolation indices $\wp_i, i = 1, 
\ldots, r$ need to be applied to the matrix $U_F$ or the nonaffine function $b(\mu)$. This implies 
that $P^TU_F$ merely consists of the rows of $U_F$ which correspond to the interpolation indices 
$\wp_i, i = 1, \ldots, r$. Similarly, $P^Tb(\mu)$ is a condensed vector composed of a few entries of 
$b(\mu)$ which correspond to the same indices. 

\Cref{myalg} provides a brief overview of the DEIM procedure.
\begin{algorithm}[H]
  \caption{DEIM algorithm \cite{morChaS10, Feng2016}}\label{myalg}
  \begin{algorithmic}[1]
      \Require POD basis $\{u_i^F\}_{i=1}^r$ for $F$ in equation (10). 
      \Ensure DEIM basis $U_F$ and indices $\vec{\wp} = [\wp_1,\ldots,\wp_r]^T \in \mathbb{R}^r$.
	\State $\wp_1$ = $\textrm{arg} \max\limits_{j \in \{1, \ldots, \mathcal{N}\}}\lvert u_{1j}^F 
		\rvert$, where $u_1^F = (u_{11}^F, \ldots, u_{1\mathcal{N}}^F)^T$. 
	\State $U_F = [u_1^F]$, $P = [e_{\wp_1}], \vec{\wp} = [\wp_1]$.
	  \For{\texttt{i = 2 to r}}
	  \State Solve $(P^TU_F)\alpha = P^Tu_i^F$ for $\alpha$, where 
	  $\alpha = (\alpha_1, \ldots, \alpha_{i-1})^T$,
	  \State $r_i = u_i^F - U_F\alpha$,
	  \State $\wp_i$ = $\textrm{arg} \max\limits_{j \in \{1, \ldots, 
		  \mathcal{N}\}}\lvert r_{ij} \rvert$, where $r_i = (r_{i1}, 
		  \ldots, r_{i\mathcal{N}})^T$. 
	  \State $U_F \gets [U_F \ u_i^F]$, $P \gets [P \ e_{\wp_i}]$, $\vec{\wp} \gets 
		\begin{bmatrix}
		\vec{\wp} \\
		\wp_i
		\end{bmatrix}$.
	  \EndFor
	\State \textbf{end for}
  \end{algorithmic}
\end{algorithm}
Note that in \Cref{myalg}, the POD basis $\{u_i^F\}_{i=1}^r$ is of great significance as an input 
basis for the DEIM procedure in two ways. First, a set of interpolation indices $\wp_i$ are 
constructed inductively based on this basis through a greedy algorithm. Secondly, an error 
analysis in \cite{morChaS10} indicates that the ordering of this basis according to the dominant 
singular values makes it the right choice for this algorithm. In step 1, the process selects the 
first interpolation index $\wp_1$ which corresponds to the location of the entry in $u_1^F$ with 
the largest magnitude. The subsequent indices in step 6, $\wp_i, i = 2, \ldots, r$, are selected 
in such a way that each of them corresponds to the location of the entry in $r$ (step 5) with the 
largest magnitude.

\subsection{DEIM Approximation Error}\label{DEIM_error}

We compute the error due to the DEIM interpolation which is to be included into the residual in the 
$\textit{a posteriori}$ error estimation. This error was first proposed in \cite{morWirSH14} for 
nonlinear dynamical systems and has also been used in \cite{Feng2016} in the context of a nonlinear 
population balance systems. We extend this idea to parametrized elliptic PDEs where the DEIM error is 
given by,
\begin{equation}
 e_{\textrm{DEIM}} = b(\mu)-\tilde{b}(\mu)= \Pi_2(I-\Pi)b(\mu),
\end{equation}
where $\Pi$ and $\Pi_2$ are oblique projectors defined as follows,
\begin{equation}\label{eq:oblique_1}
 \Pi = U_F(P^TU_F)^{-1}P^T,
\end{equation}
and 
\begin{equation}\label{eq:oblique_2}
 \Pi_2 = (I-\Pi)\tilde{U}_F(\tilde{P}^T(I-\Pi)\tilde{U}_F)^{-1}\tilde{P}^T.
\end{equation}
In equation (\ref{eq:oblique_1}), $U_F = (u_1^F, \ldots, u_r^F) \in \mathbb{R}^{\mathcal{N}\times r}$ 
and $P \in \mathbb{R}^{\mathcal{N}\times r}$ are the current DEIM basis and interpolation index matrix
obtained from \Cref{myalg}. 

To obtain $\Pi_2$ in (\ref{eq:oblique_2}), we assume that $r^*(\geq r)$ DEIM basis vectors 
$U_F^* = (u_1^F, \ldots, u_{r^*}^F)$ interpolate $b(\mu)$ exactly, i.e.
\begin{equation}
 b(\mu) = U_F^*((P^*)^TU_F^*)^{-1}(P^*)^Tb(\mu),
\end{equation}
where $P^*$ is the corresponding index matrix with $r^*$ columns. Finally, 
$\tilde{U}_F = U_F^*(:,r+1:r^*)$ and $\tilde{P} = P^*(:,r+1:r^*)$ such that $U_F^* = [U_F,\tilde{U}_F]$ 
and $P^* = [P,\tilde{P}]$, where $M(:,r+1:r^*)$, using MATLAB notation \cite{Feng2016}.
In the next subsection, we introduce an \textit{a posteriori} error estimation derived from the 
residual of the approximate RB solution and the DEIM approximation error.

\subsection{\textit{A Posteriori} Error Estimation}\label{Error_estmn}

\textit{A posteriori} error estimators are computable indicators which provide an estimate to the 
actual solution error by utilizing the residual of the approximate RB solution. An efficient 
error estimator is required to possess three major characteristics, namely: it is required to be as 
sharp as possible (close to the unknown actual error), asymptotically correct (tend to zero with 
increasing RB space dimension $N$, at a similar rate as the actual error), and computationally cheap 
(because it is computed in the online phase). Therefore, these estimators guarantee both reliability 
and efficiency of the reduction process \cite{morQuart16}. 

We first compute the residual due to DEIM interpolation; 
\begin{equation}\label{eq:DEIM_residual}
 r_N^{\textrm{DEIM}}(u_N;\mu) = (\rho + \tilde{b}(\mu)) - A^{\mathcal{N}}(\mu)u_N(\mu),
\end{equation}
where $\tilde{b}(\mu) = \Pi b(\mu)$ is the DEIM interpolation of $b(\mu)$ and 
$u_N(\mu):= Vu^N(\mu)$ is the RB solution transformed back to the high-fidelity space 
$\mathcal{N}$. Then the final residual is obtained by including the DEIM approximation error 
derived in \Cref{DEIM_error} as follows;
\begin{equation} \label{eq:General_residual}
   \begin{split}
 r_N(u_N;\mu) & = (\rho + b(\mu)) - A^{\mathcal{N}}(\mu)u_N(\mu) \\
	     & = (\rho + \tilde{b}(\mu)) - A^{\mathcal{N}}(\mu)u_N(\mu) + b(\mu) - \tilde{b}(\mu)\\
	     & = r_N^{\textrm{DEIM}}(u_N;\mu) + \underbrace{b(\mu) - \tilde{b}(\mu)}_{:=e_{\textrm{DEIM}}}\\
	     & = r_N^{\textrm{DEIM}}(u_N;\mu) + e_{\textrm{DEIM}}.
   \end{split}
\end{equation}
The \textit{a posteriori} error estimation is then derived from the residual in 
(\ref{eq:General_residual}). Rewriting the first equation of (\ref{eq:General_residual}), we 
obtain
\begin{equation} \label{eq:residual}
  \begin{split}
  r_N(u_N;\mu) & = A^{\mathcal{N}}(\mu)u^{\mathcal{N}}(\mu) - A^{\mathcal{N}}(\mu)u_N(\mu) \\
                & = A^{\mathcal{N}}(\mu)e(\mu), 
  \end{split}
\end{equation}
where the error $e(\mu):= u^{\mathcal{N}}(\mu) - u_N(\mu)$ is given by
\begin{equation}\label{eq:error}
 e(\mu)  = (A^{\mathcal{N}}(\mu))^{-1}r_N(u_N;\mu).
\end{equation}
We obtain an upper bound for the 2-norm of the error by taking the 2-norm on both sides of 
equation (\ref{eq:error}), i.e.
\begin{equation}\label{eq:error_norm}
\begin{split}
 \Arrowvert e(\mu)\Arrowvert_2 \leq \Arrowvert (A^{\mathcal{N}})^{-1}(\mu)\Arrowvert_2
 \Arrowvert r_N(u_N;\mu)\Arrowvert_2 &= \frac{\Arrowvert r_N(u_N;\mu)\Arrowvert_2}{\sigma_{min}(A^{\mathcal{N}}(\mu))}
  \\& =: \tilde \Delta_N(\mu),
  \end{split}
\end{equation}
where $\sigma_{min}(A^{\mathcal{N}}(\mu))$ is the smallest singular value of 
$A^{\mathcal{N}}(\mu)$ \cite{morQuart16}. 
The quantity $\tilde \Delta_N(\mu)$ is a rigorous error bound, and can be used to select 
snapshots within the greedy algorithm in the offline stage and consequently to measure the 
accuracy of the RB approximation in the online stage \cite{morHesB13}. For efficient computation of 
the norm of the residual and error bounds, see \cite{morHest16,morQuart16}.  
It is computationally expensive to compute $\sigma_{min}(A^{\mathcal{N}}(\mu))$ in the online 
phase as it entails solutions of large-scale eigenvalue problems \cite{morHesB13}. Therefore, in 
our computations, we use the norm of the residual as our error estimator, which satisfies the 
inequality (\ref{eq:error_norm}) and provides an estimation of the true error that works well for 
our problem. It also provides rapid convergence as depicted in the numerical results in 
\Cref{fig:True_Est}. It is given by
\begin{equation}
 \Arrowvert e(\mu)\Arrowvert_2 \approx \Arrowvert r_N(u_N;\mu)\Arrowvert_2 := \Delta_N(\mu).
\end{equation}

\section{Numerical Results}\label{Numerical_results}

\subsection{Finite Difference Results}

We consider the LPBE (\ref{eq:LPBE}), a parameter domain $\mu \in D = [0.05,0.15]$, and a cubic grid 
of $129$ points and a box length of $60\,\textrm{\AA}$ centered at the protein position. The 
parameter domain is chosen for a feasible physiological process and $\mu$ resides in the second 
term in the kappa function. Information about the molecular charge density is obtained from a PQR 
file which contains $1228$ atoms of the protein \textit{fasciculin 1} toxin CPDB entry 1FAS. We discretize 
the LPBE with a centered finite difference scheme and the resulting parametrized linear system 
(\ref{eq:FOM}) is of more than $2\times 10^{6}$ degrees of freedom. This FOM is solved by the 
aggregation-based algebraic multigrid (AGMG) method, where a tolerance of $10^{-10}$ and 
a zero initial guess are used \cite{Notay:10,NapNot:12,Notay:12}. 

The choice of the tolerance directly affects the results of the greedy algorithm. 
Therefore, it is prudent to ensure that the high-fidelity solution ($u^\mathcal{N}(\mu), 
\mathcal{N}= 2,146,689$) is highly accurate. Some of the iterative methods commonly used in the PBE solvers 
are; the minimal residual (MINRES) method, the generalized minimal residual (GMRES) method and the 
biconjugate gradient stabilized (BICGSTAB) method. These methods employ the incomplete LU factorization 
to generate the preconditioner matrices L (lower diagonal) and U (upper diagonal) which are used to improve 
their stability and convergence at low costs \cite{Vergara-Perez2016}.

\begin{figure}[t]
\begin{tabular}{@{}cc@{}}
\includegraphics[height=1.25in]{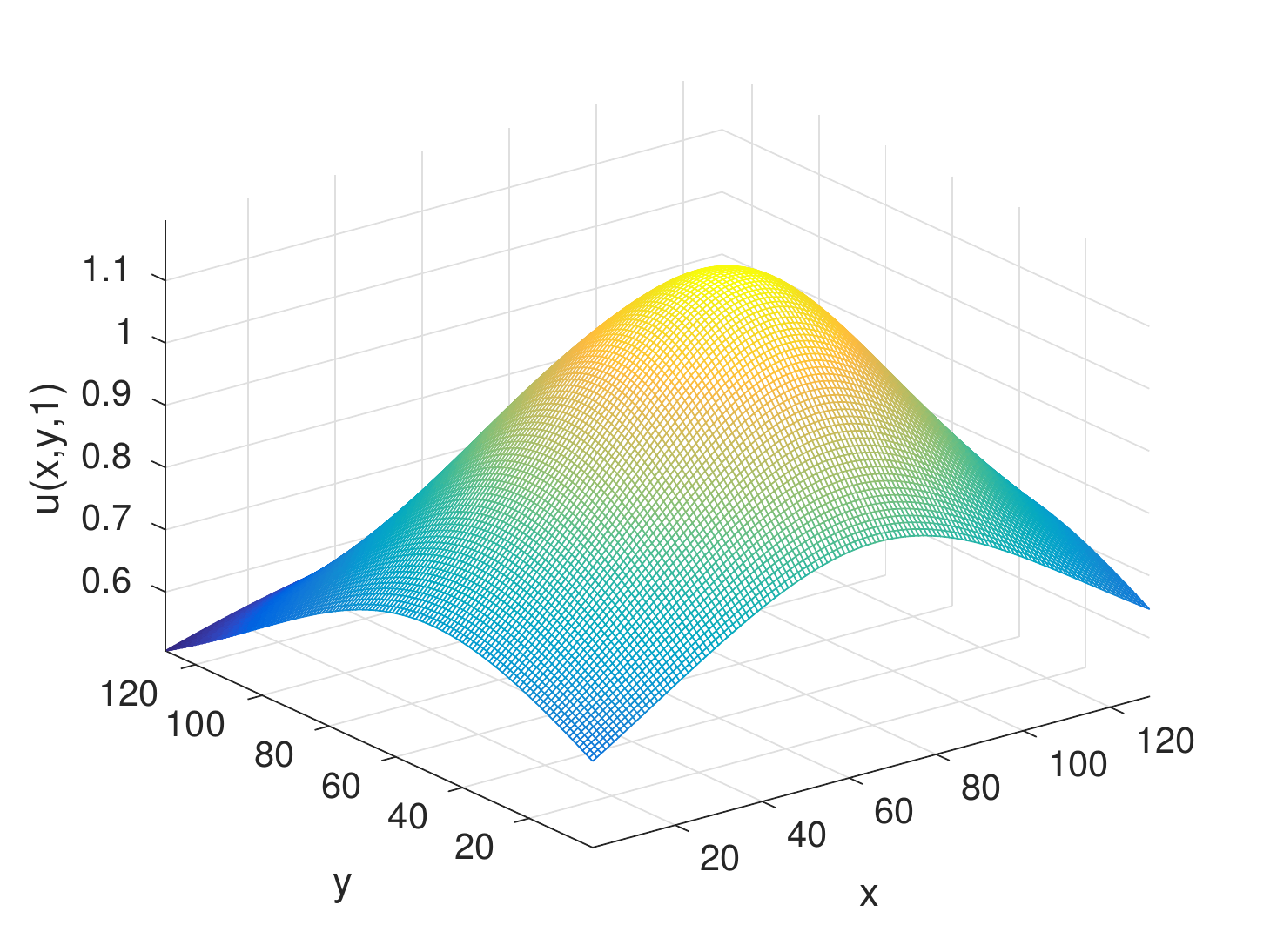} & \includegraphics[height=1.25in]{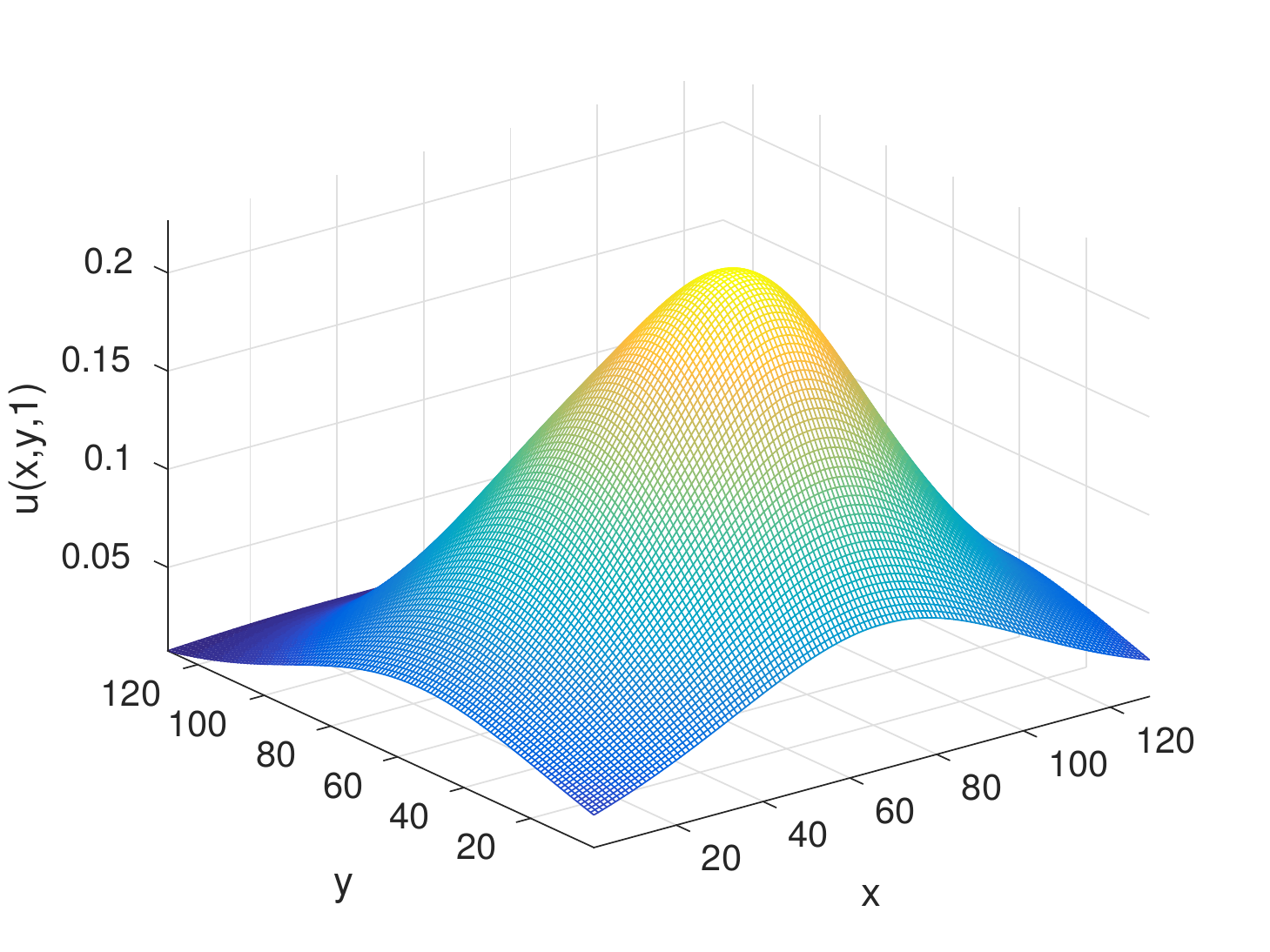}\\
\includegraphics[height=1.25in]{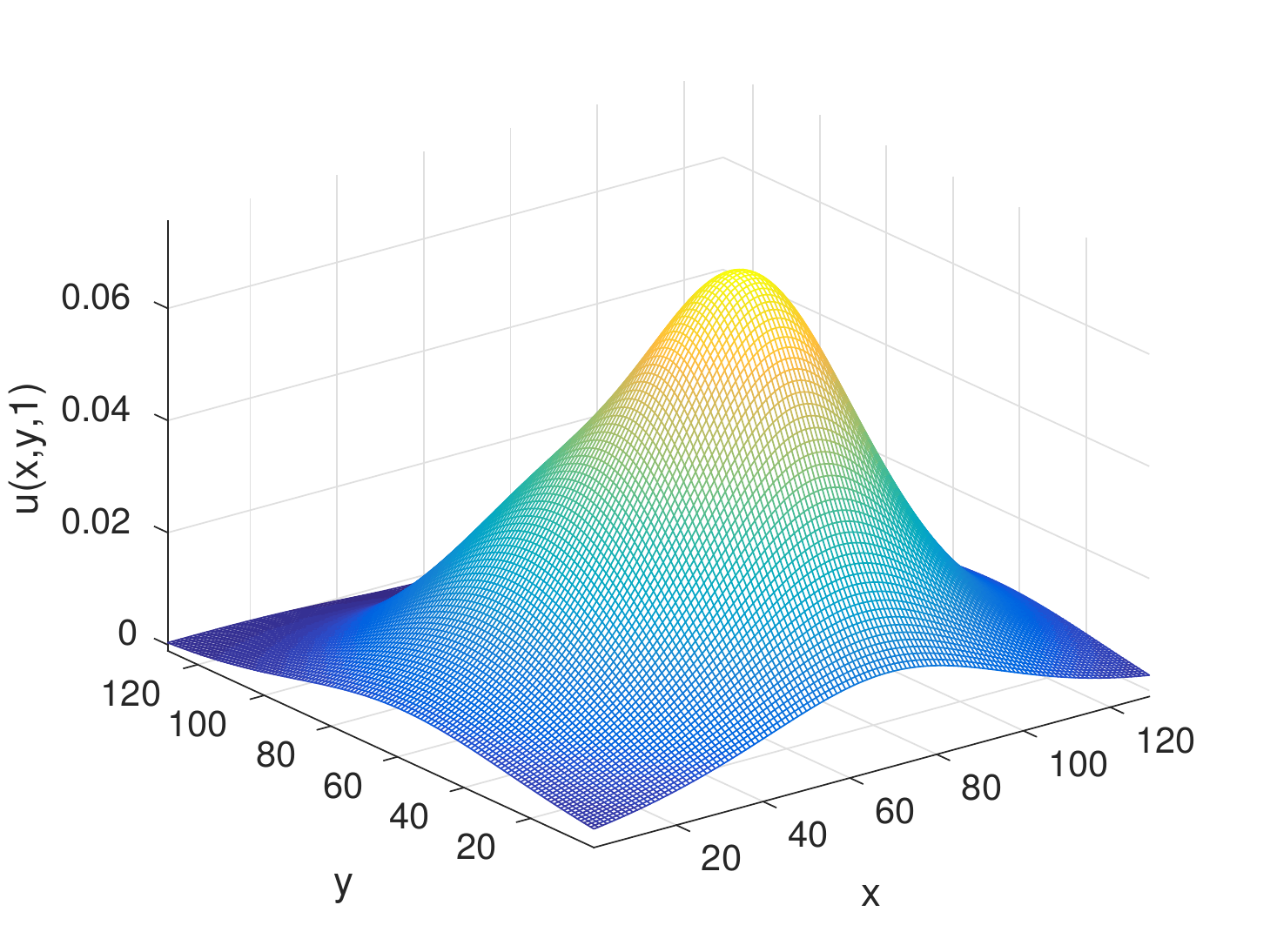} & \includegraphics[height=1.25in]{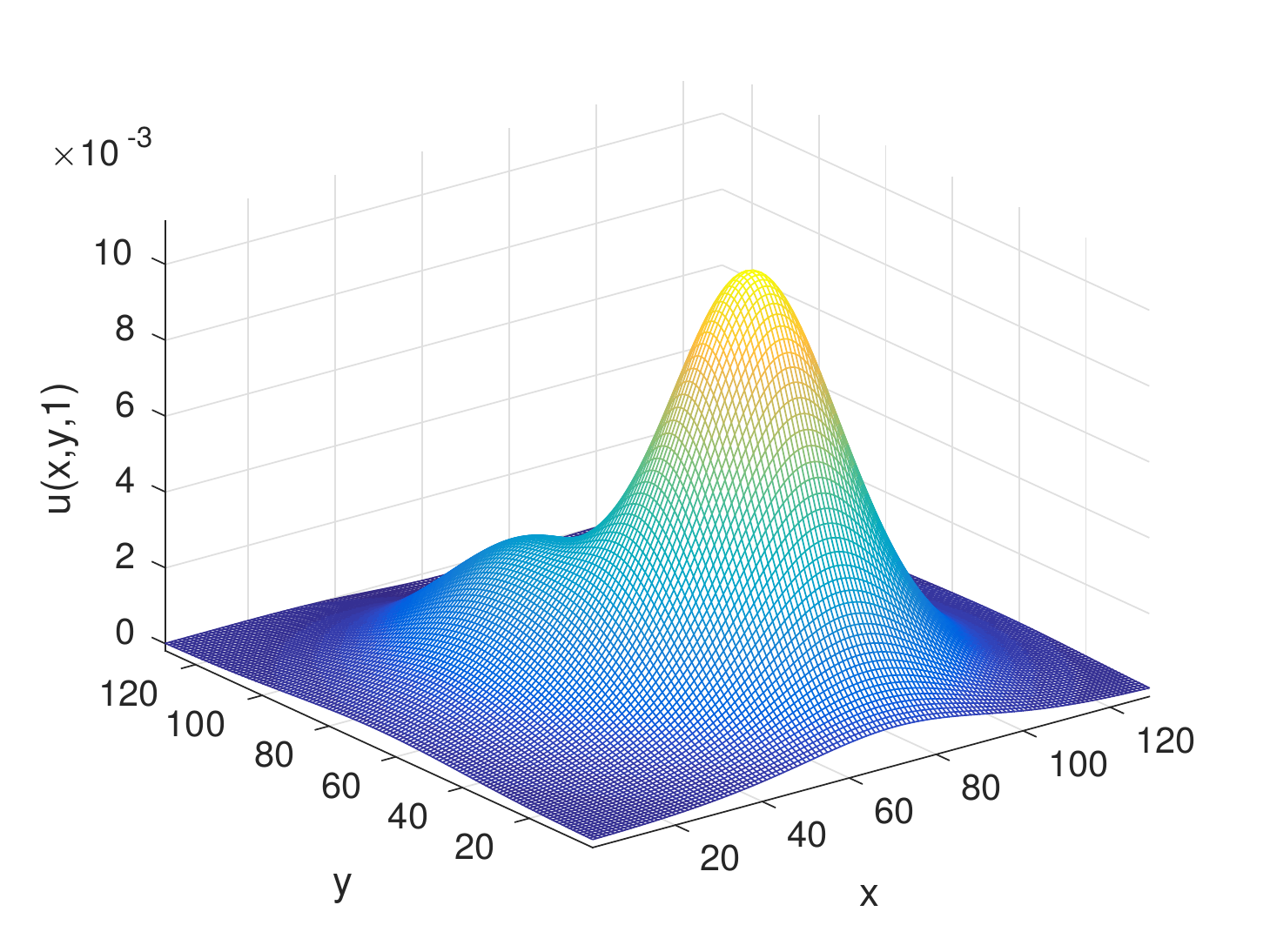}
\end{tabular}
\caption{\label{fig:NumsolsB}%
High-fidelity solutions ($u^{\mathcal{N}}(\mu)$) at varying ionic \\strengths (i.e., 
$\mu$ = \{0, 0.05, 0.15, 0.5\}), respectivley.}
\end{figure}

\Cref{fig:NumsolsB} shows the lower cross-sections of the $z$-axis 
of the electrostatic potential $u(x,y,1)$. Note that the potential decays exponentially with the
variation of the parameter $\mu$, and is attributed to the large force constant ($332$ kcal/mol) 
of electrostatic interactions. In the absence of ions (that is, at $\mu = I = 0$), these 
interactions are long ranged, but in the presence of ions (that is, $\mu > 0$), they are damped 
or screened and gradually decay to zero \cite{Fogolari2002}. The computational time taken to 
obtain the high-fidelity solution $u^\mathcal{N}(\mu)$ is approximately $28$ seconds on average 
and varies depending on the value of the ionic strength used.

\subsection{Accuracy of FDM} \label{FDM_accuracy}

We demonstrate the accuracy and reliability of the FDM before applying the RBM for the solution 
of the PBE. This is because the accuracy of the RBM depends on that of the underlying discretization 
technique. In this study, we consider six test examples to validate the FDM which include a Born 
ion and five proteins consisting of between $380$ and $3400$ atoms, respectively. We compare the FDM 
results with those of APBS for electrostatic solvation free energy at different mesh refinements. Firstly, 
we consider the Born ion which is a canonical example for polar solvation and whose analytical 
solution is well known. 

This analytical solution gives the polar solvation energy which results from the transfer of a 
non-polarizable ion between two dielectrics \cite{www}, i.e.,
\begin{equation}\label{eq:Born_energy}
 \Delta_pG_{Born} = \frac{q^2}{8\pi \epsilon_0 r}(\frac{1}{\epsilon_{out}}-\frac{1}{\epsilon_{in}}),
\end{equation}
where $q$ is the ion charge, $r$ is the ion radius, $\epsilon_{out}$ is the external dielectric 
coefficient (e.g., water) and $\epsilon_{in}$ is the internal dielectric coefficient (e.g., vacuum). 
This model assumes zero ionic strength. We consider a Born ion of unit charge, $3\textrm{\AA}$ radius 
and located at the origin ($(0,0,0)$). Here, $\epsilon_{in} = 1$ and $\epsilon_{out} = 78.54$. With 
these parameters, the analytical solution in (\ref{eq:Born_energy}) is
\begin{equation}
\Delta_pG_{Born} = -691.85(\frac{q^2}{r})= -230.62 \textrm{kJ/mol}.
\end{equation}

We compare numerical computations using equation (\ref{eq:Elec_energy}) for charging free energies in a
homogeneous ($\epsilon_{in} = \epsilon_{out} = 1$) and heterogeneous ($\epsilon_{in} = 1, 
\epsilon_{out} = 78.54$) dielectric coefficients with the analytical solution \cite{www}. We use the 
following additional parameters. We consider two different mesh sizes (or $\Delta x$), which result in different 
degrees of freedom (or $\mathcal{N}$) as shown in \Cref{table:Energy_compare_born}. Numerical results using 
FDM are compared with those of the exact solution (\ref{eq:Born_energy}) and APBS (which uses FEM). The 
results show that the FDM method gives solutions which are consistent with those of the exact solutions, 
as well as those of the APBS software package.
\begin{table}[H]
\centering
\begin{tabular}{|c|c|c|c|c|c|}
  \hline
  $\Delta x$ & $\mathcal{N}$ & Solver & Numerical 
  & Analytical & Relative error\\ \hline
  \multirow{2}{*}{0.33} & \multirow{2}{*}{$97^3$} & APBS & -229.59 & -230.62 & 4.4662e-3\\ 
   & ~ & FDM & -232.86 & -230.62 & 9.7130e-3\\ \hline
  \multirow{2}{*}{0.25} & \multirow{2}{*}{$129^3$} & APBS & -230.00 & -230.62 & 2.6884e-3\\ 
   & ~ & FDM & -230.42 & -230.62 & 8.6723e-4\\ \hline
\end{tabular}
\caption{Comparison of Born ion solvation energies in kJ/mol.}
\label{table:Energy_compare_born}
\end{table}

Secondly, we compare the accuracy of FDM for the LPBE with the following set of typical examples of use of 
LPBE and APBS in particular: Calculation of the total electrostatic energy (including self-interaction 
energies) of a $22$ residue, $\alpha$-helical peptide from the N protein of phage $\lambda$ which binds 
to its cognate $19$ nucleotide box B RNA hairpin \cite{GGDr:03}, Fasciculin 1, an anti-acetylcholinesterase 
toxin from green mamba snake venom \cite{DuMaBoFo:92}, the electrostatic potential of a minimized FKBP 
protein from binding energy calculations of small ligands \cite{BurkTayWalk:00}, a 180-residue cytokine 
solution NMR structure of a murine-human chimera of leukemia inhibitory factor (LIF) 
\cite{HinMauZhaNic:98}, and the binding energy of a balanol ligand to the catalytic subunit of the 
CAMP-dependent protein kinase A, here the apo form of the enzyme \cite{Narayana:99}. The proteins and 
or complexes have the following number of atoms (379, 1228, 1663, 2809, and 3423), respectively. 

The electrostatic solvation free energies, $\Delta E$ are computed and shown in 
\Cref{table:Energy_compare} for varying grid resolutions $\Delta x$. 
However, we here do not have the analytical electrostatic energies for these proteins but rely on the 
accuracy of the APBS software for validation. A compute cluster with 4 Intel Xeon E7-8837 CPUs running 
at 2.67 GHz (8 cores per CPU) and 1 TB RAM, split into four 256 GB parts (each CPU controls one part) 
is used to carry out the computations which require a huge amount of memory, so that it allows for 
solving large-scale problems with $\mathcal{N} \geq (3\times 10^6)$.
\begin{table}[t]
\centering
\begin{tabular}{|K{1.4cm}|K{1.3cm}|K{1.7cm}|K{1.7cm}|K{1.6cm}|}
  \hline
  $\Delta x$& $\mathcal{N}$ & $\Delta E$, FDM & $\Delta E$, APBS 
  & Relative error\\ \hline
   \multicolumn{5}{|p{\linewidth}|}{1. Solvation energies of a $22$ residue, $\alpha$-helical 
   peptide from the N protein of phage $\lambda$ in kJ/mol. (379 atoms)} \\ \hline
  {0.375} & $129^3$ & -4557.7052 & -4546.5150 & 2.4613e-3\\ \hline
  {0.320} & $161^3$ & -4541.4782 & -4532.7595 & 1.9235e-3\\ \hline
  {0.260} & $193^3$ & -4522.4752 & -4516.8544 & 1.2444e-3\\ \hline
  \multicolumn{5}{|p{\linewidth}|}{2. Solvation energies of \textit{fasciculin 1} in kJ/mol. 
  (1228 atoms)} \\ \hline
  {0.465} & $129^3$ & -5870.5357 & -5845.8594 & 4.2212e-3\\ \hline
  {0.375} & $161^3$ & -5684.8448 & -5664.8475 & 3.5301e-3\\ \hline
  {0.320} & $193^3$ & -5629.1979 & -5611.2503 & 3.1985e-3\\ \hline
  \multicolumn{5}{|p{\linewidth}|}{3. Solvation energies of the electrostatic potential of a 
  minimized FKBP protein in kJ/mol. (1663 atoms)} \\ \hline
  {0.465} & $129^3$ & -4419.0384 & -4403.8761 & 3.4429e-3\\ \hline
  {0.375} & $193^3$ & -4344.5491 & -4331.1010 & 3.1050e-3\\ \hline
  {0.320} & $225^3$ & -4292.5359 & -4288.0842 & 1.0382e-3\\ \hline
  \multicolumn{5}{|p{\linewidth}|}{4. Solvation energies of a 180-residue cytokine solution NMR 
  structure of a murine-human chimera of leukemia inhibitory factor (LIF) in kJ/mol. (2809 atoms)} \\ \hline
  {0.450} & $161^3$ & -9317.7636 & -9293.9750 & 2.5595e-3\\ \hline
  {0.375} & $193^3$ & -9270.0472 & -9247.2822 & 2.4618e-3\\ \hline
  {0.280} & $257^3$ & -9153.9477 & -9134.2879 & 2.1523e-3\\ \hline
  \multicolumn{5}{|p{\linewidth}|}{5. Solvation energies of CAMP-dependent protein kinase A, 
  here the apo form, in kJ/mol. (3423 atoms)} \\ \hline
  {0.465} & $129^3$ & -19742.3639 & -19681.3183 & 3.1017e-3\\ \hline
  {0.375} & $161^3$ & -19332.6588 & -19296.6336 & 1.8670e-3\\ \hline
  {0.320} & $193^3$ & -19039.8581 & -19014.0380 & 1.3579e-3\\ \hline
\end{tabular}
\caption{Comparison of electrostatic solvation free energies \\$\Delta E$, between FDM and APBS for 
different proteins.}
\label{table:Energy_compare}
\end{table}

From Table \ref{table:Energy_compare}, we can clearly see that the 
results of the FDM method agree well with those of APBS in terms of convergence with respect to mesh 
refinement. Hence, we conclude that we can test the RBM in conjunction with our FDM solver reliably. We 
expect no differences when using a FEM solver like APBS, which would require intruding the software.

\subsection{Accuracy of the RBM}

\begin{figure}[t]
\centering
\begin{tikzpicture}
\begin{customlegend}[legend columns=-1, legend style={/tikz/every even column/.append 
style={column sep=2.0cm}} , legend entries={True error, Error estimator}, ]
\addlegendimage{blue,dashed, line width = 2pt}
\addlegendimage{black,solid, line width = 2pt}
\end{customlegend}
\end{tikzpicture}

\begin{tabular}{@{}cc@{}}
\includegraphics[height=1.25in]{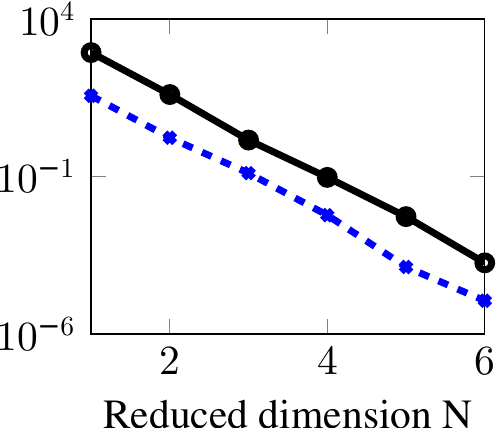} & 
\includegraphics[height=1.25in]{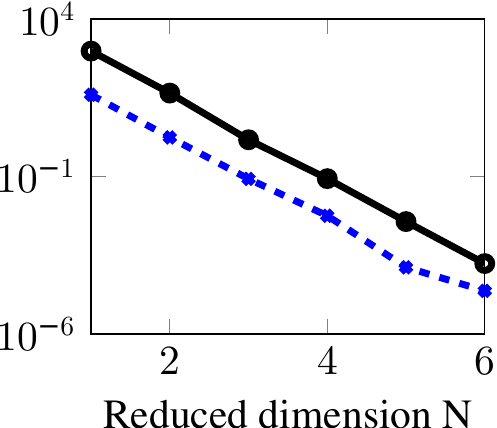}\\
\includegraphics[height=1.25in]{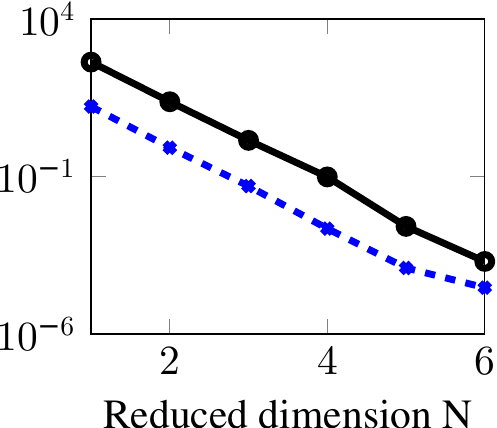} & 
\includegraphics[height=1.25in]{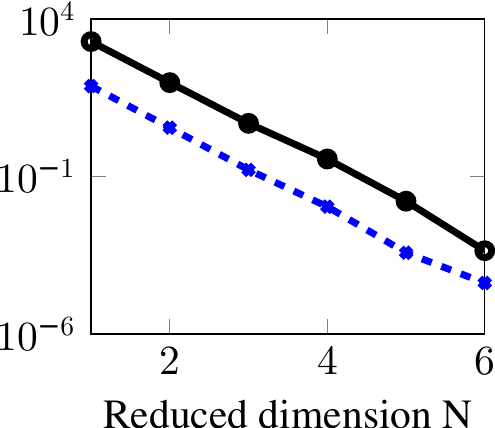}\\
\includegraphics[height=1.25in]{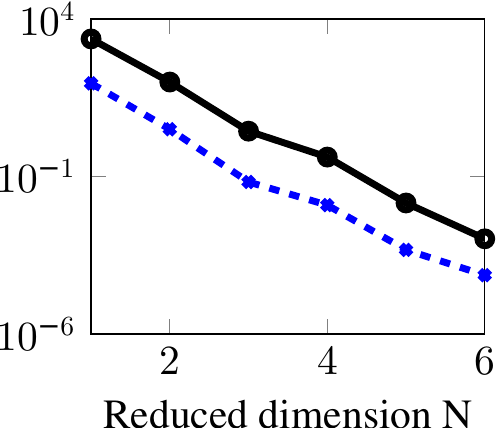}
\end{tabular}
\caption{Comparison of maximal error estimator and \\true error for the proteins in  
\Cref{table:Energy_compare}, respectively.}
\label{fig:True_Est}
\end{figure}

In this section, we evaluate the accuracy of the RBM for the approximation of the high-fidelity 
solutions generated by the FDM for the five proteins which were investigated in \Cref{FDM_accuracy}. We 
consider a cubic domain of $129$ points and a box length of $60\,\textrm{\AA}$ centered at the protein 
position for all the computations. \Cref{fig:True_Est} shows the decay 
of the error estimator and the true error during the greedy algorithm at the current RB dimension 
$i = 1, \ldots, N$. They corroborate the asymptotic correctness property stated in \Cref{Error_estmn}, 
and it is evident that the error estimator is an upper bound to the true error. We also observe a high 
convergence rate of the error estimator with up to two orders of magnitude and the RB space is rich 
enough at only six iterations of the greedy algorithm for the five proteins. These error estimators are 
the maximal error and relative maximal error, respectively, and are defined as, 
$\Delta_N^{\max} = \max\limits_{\mu\in \Xi}\Arrowvert r_N(u_N;\mu)\Arrowvert_2,$ and 
${\Delta_N^{\max}}/{\Arrowvert u_N(\mu^*)\Arrowvert_2},$
where $\mu^*=\arg\max\limits_{\mu\in \Xi}\|r_N (u_N; \mu)\|_2.$ 

In the greedy algorithm, we apply an error tolerance of $\epsilon = 10^{-3}$ and a training set $\Xi$ 
consisting of $l=11$ samples of the parameter. From \Cref{fig:True_Est}, it is evident that both the 
error estimator and the true error fall below the prescribed tolerances at the final dimension of the 
ROM (i.e. $N = 6$). 

\begin{figure}[t]
\centering
\begin{tikzpicture}
    \begin{customlegend}[legend columns=-1, legend style={/tikz/every even column/.append 
    style={column sep=2.0cm}} , legend entries={True error, Error estimator}, ]
    \addlegendimage{blue,solid, line width = 2pt}
    \addlegendimage{black,solid, line width = 2pt}
    \end{customlegend}
  \end{tikzpicture}
 
 \begin{tabular}{@{}cc@{}}
\includegraphics[height=1.25in]{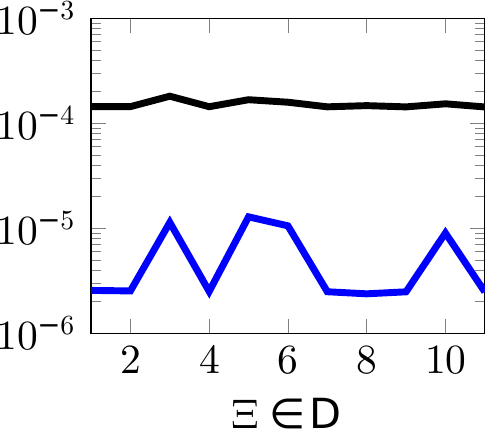} & 
\includegraphics[height=1.25in]{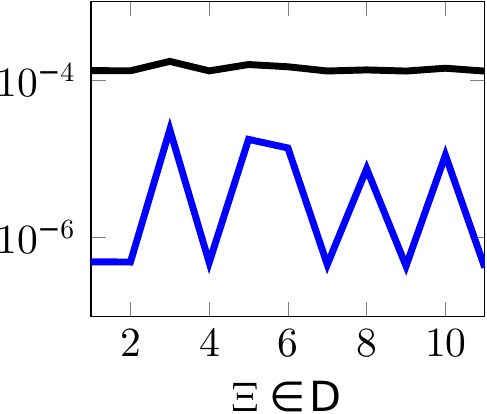}\\
\includegraphics[height=1.25in]{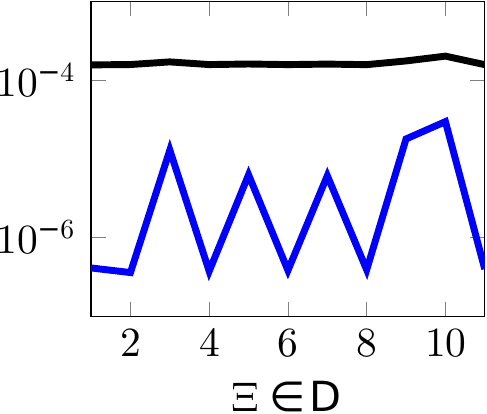} & 
\includegraphics[height=1.25in]{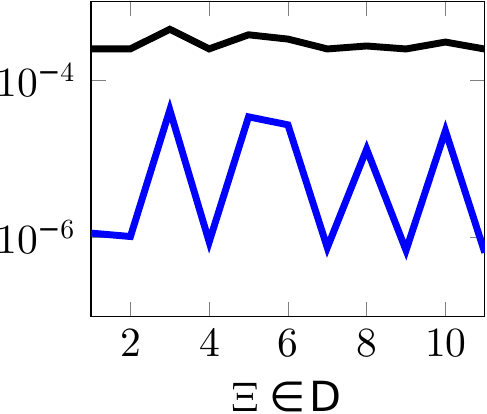}\\
\includegraphics[height=1.25in]{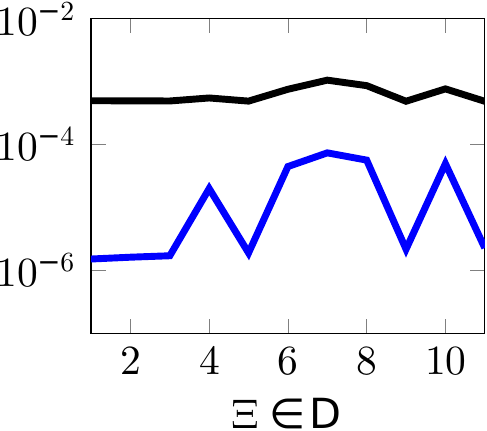}
\end{tabular}
\caption{Comparison of error estimator and true error for \\the final ROM for $\Xi \in D$ and for the 
proteins in \Cref{table:Energy_compare}, \\respectively.}
\label{fig:Truedelt_err}
\end{figure}

\Cref{fig:Truedelt_err} shows the error estimator and the true error of the finally constructed ROM 
over $\mu_i = \Xi$, for $i = 1, ..., 11$ samples for each protein as in \Cref{table:Energy_compare}, 
respectively. It is evident that the error estimator for the final RB approximations of dimension $N = 6$ 
is indeed an upper bound of the true error and a trend that both quantities behave similarly is clearly 
visible from the graphs. Consequently, the error estimators fall below the greedy tolerance of 
$10^{-3}$.

\begin{figure}[t]
\centering
\begin{tikzpicture}
    \begin{customlegend}[legend columns=-1, legend style={/tikz/every even column/.append 
    style={column sep=2.0cm}} , legend entries={True error}, ]
    \addlegendimage{blue,solid, line width = 2pt}
    \end{customlegend}
  \end{tikzpicture}
  
  \begin{tabular}{@{}cc@{}}
\includegraphics[height=1.25in]{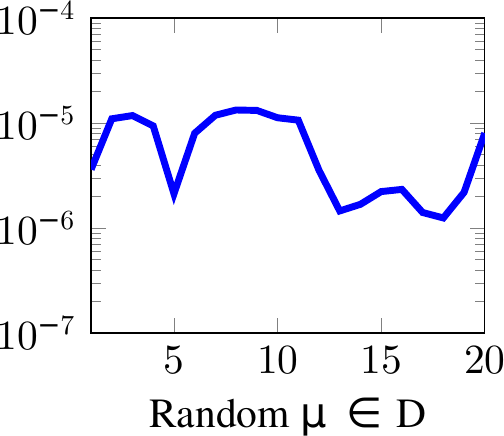} & 
\includegraphics[height=1.25in]{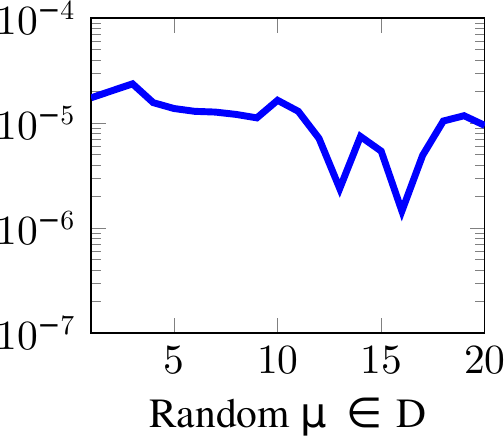}\\
\includegraphics[height=1.25in]{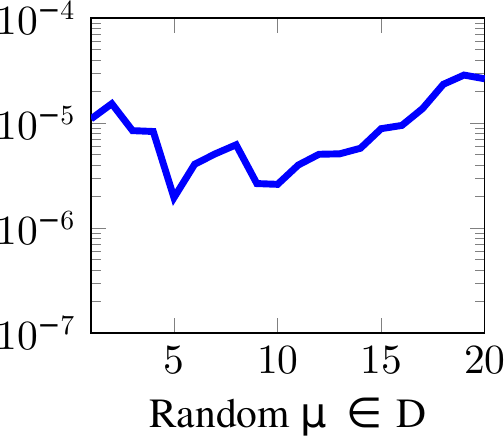} & 
\includegraphics[height=1.25in]{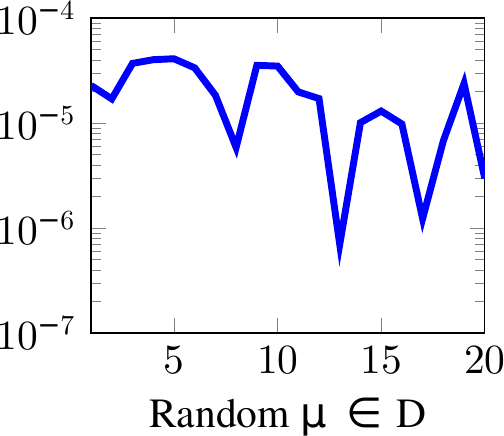}\\
\includegraphics[height=1.25in]{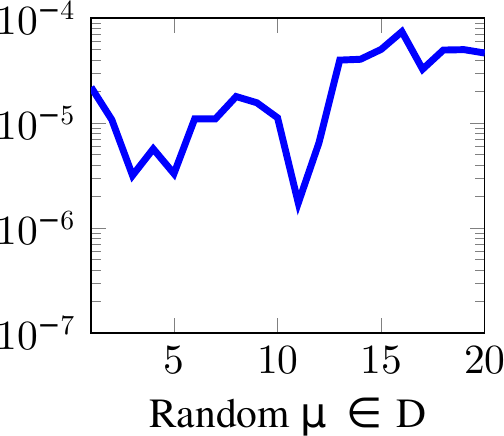}
\end{tabular}  
\caption{True error $\|u^{\mathcal{N}}(\mu)-u_N(\mu)\|_2$ for random parameters\\ $\mu \in D$ for the proteins 
in \Cref{table:Energy_compare}, respectively.}
\label{fig:Truedelt_err_rand}
\end{figure}

\Cref{fig:Truedelt_err_rand} is used to validate the true error in \Cref{fig:Truedelt_err}, 
whereby $20$ random values of the parameter domain $D$ which are different from those in the training set 
$\Xi$ are used. A common observation from these figures is that the true errors fall below  
$\mathcal{O}(10^{-4})$, which is an order of magnitude below the error estimator. The computational 
time taken to obtain the approximate solution $u_N(\mu)$ in the online phase is approximately 
$4.97\times 10^{-3}$ seconds on average, for any parameter $\mu \in D$.

\subsection{Runtimes and Computational Speed-ups}

Before we dive into the runtimes of the various phases of the RBM, we would like to make clear about some key 
notions of the two phases of the greedy algorithm, i.e., the offline and online phases, respectively. The 
offline phase is subdivided into two parts, the offline-offline phase, and the offline-online phase 
\cite{morHest16}. The offline-offline phase involves computation of the snapshots and pre-computing the 
parameter-independent quantities. The offline-online phase involves computation of the error estimator and 
the RB approximation. On the other hand, the pure online phase is where the final ROM has been constructed 
after the accuracy of the reduced basis is fulfilled, and is independent of the greedy algorithm. In this 
phase, the ROM can be solved for any parameter value in the parameter domain, including those which are 
different from the training set.

Table \ref{table:speed_up} shows the runtimes and computational speed-ups obtained with the use of DEIM 
approximation during the offline-online phase of the RBM at a single iteration of the greedy algorithm 
and with the use of the RBM in solving the linear system. We use a modest PC with 
the following specifications: Intel (R) Core (TM)2 Duo CPU E8400 @ 3.00GHz with 8GB RAM. In this section, 
the PBE is applied to the protein \textit{fasciculin 1}.
\begin{table}[H]
\centering
\begin{tabular}{|c|c|c|c|}
  \hline
  \multicolumn{4}{|c|}{Runtime (seconds) and speed-up} \\
  \hline
   & Without DEIM & With DEIM & Speed-up\\
  \hline
  Offline-online phase & 96.29 & 4.84 & 20\\ \hline
  Assemble and solve & & &\\ ROM & 8.36 & 9.91e-03 & 844\\ \hline
\end{tabular}
\caption{Runtimes and speed-ups due to DEIM.}
\label{table:speed_up}
\end{table}

Table \ref{table:speed_up2} shows the runtimes of computing the FOM and the ROM at a given parameter value, 
respectively. The runtimes at different phases of the RBM are also presented. Speed-up factors induced by 
solving the ROM are listed to visualize the big difference between the FOM and the ROM. The ROM is much 
faster and takes a split second to assemble and solve for any parameter value. In the offline phase of the 
RBM, which comprises the greedy algorithm, the dominating cost is that of solving the linear system of 
the FOM by AGMG (i.e., computing snapshot) at every iteration of the greedy algorithm. Miscellaneous in 
this case refers to the runtime to initialize the FDM, including assembling the FOM. The total RBM 
runtime includes the miscellaneous and offline runtimes. 
\begin{table}[H]
\centering
\begin{tabular}{|c|c|c|c|}
  \hline
  \multicolumn{4}{|c|}{Runtime (seconds) and speed-up} \\
  \hline
   & FOM & ROM & Speed-up\\
  \hline
  Solve linear system & 11.88 & 4.97e-03 & 7,616\\ \hline
  Assemble and solve linear system & 27.82 & 9.91e-03 & 5,500\\ \hline
   \multicolumn{4}{|c|}{Runtime (seconds) for RBM phases} \\
  \hline
   Miscellaneous & Offline & Online & Total RBM\\ \hline
    10.58  & 85.54 & 9.91e-03 & 96.12 \\ \hline
\end{tabular}
\caption{Runtimes and speed-ups for FOM, ROM and RBM.}
\label{table:speed_up2}
\end{table}

Table \ref{table:runtimes} shows the runtimes of APBS and RBM for solving the FOM and the ROM at any given 
parameter value, respectively. The speed-up factor of RBM w.r.t. the APBS is also shown for different 
numbers of parameter values. It is evident that RBM is much more efficient than APBS when solving the 
system for many input parameter values (i.e. in a multi-query context). This is because we only need to 
solve a small system of order $N = 6$ once the final ROM model has been constructed which takes 
approximately $9.91\times 10^{-3}$ seconds for each parameter value, whereas APBS solves the FOM besides 
the initial system setup. 

In a nutshell, to solve the LPBE for any parameter value with APBS, it takes 
$22.893$ seconds, because the solver has to reconstruct the linear system. This implies that it takes 
approximately $2,289.3$ seconds to compute the potential for $100$ parameter values (neglecting the 
runtime to modify the input files). This is more expensive than the total RBM time of $96.12$ seconds. 
On the other hand, it takes the RBM approximately $9.91\times 10^{-1}$ seconds to solve the ROM of the 
LPBE for the same number of parameters values (i.e., 100). 

The RBM only solves the FOM $N$ times during the expensive offline phase as stated in \Cref{euclid}. Moreover, 
the RBM utilizes the precomputed system matrices and vectors and only solves the ROM for the new parameter 
value, thus saving a significant amount of computational costs during the online phase. This efficient 
implementation of a new mathematical approach to solve the PBE holds great promise towards reducing 
computational costs in a multi-query scenario and molecular dynamics simulation. 
\begin{table}[H]
\centering
\begin{tabular}{|c|c|c|c|}
  \hline
  \multicolumn{4}{|c|}{Runtime (seconds) and speed-up for APBS and RBM} \\
  \hline
  No. of parameters & APBS & RBM & Speed-up \\
  \hline
  1 & 22.893 & $\approx 96.12$ & 0.24 \\ \hline
  10 & 228.93 & $\approx 96.12$ & 2.38 \\ \hline
  100 & 2,289.3 & $\approx 96.12$ & 24 \\ \hline
  1000 & 22,893 & $\approx 106.12$ & 215.75 \\ \hline
\end{tabular}
\caption{Runtimes for APBS and RBM.}
\label{table:runtimes}
\end{table}

\section{Conclusions}

In this paper, we have presented a new, computationally efficient approach to solving the PBE for varying 
parameter values. The RBM reduces the high-dimensional full order model by a factor of 
approximately $360,000$ and the computational time by a factor of approximately $7,600$. The error 
estimator provides fast convergence to the reduced basis approximation at an accuracy of 
$\mathcal{O}(10^{-3})$. 
The true error between the RBM and the FDM is smaller than $O(10^{-4})$, for all the parameter samples 
tested. DEIM provides a speed-up of 20 in the online phase by reducing the complexity of the nonaffine 
Dirichlet boundary conditions. This is achieved by only selecting a few entries from a high-dimensional 
vector which provides the most important information. Therefore, the RBM can be extremely beneficial in 
cases where simulations of the PBE for many input parameter values are required. This method can also be 
implemented in the available PBE solvers, for example, APBS, after a few adjustments regarding 
parametrization in the linear system are made. Our future research is based on two aspects. Firstly, we 
plan to develop a more efficient error estimator which is more rigorous than merely taking the norm of 
the residual. Secondly, we aim to develop a modified version of the LPBE which considers the PBE as 
interface problem by applying a range-separated tensor format. This is expected to reduce the 
computational complexity experienced by the current PBE studies, and to provide more accurate results 
due to the more realistic model.

\begin{acknowledgments}
The authors thank the following organizations for financial and material support on 
this project: International Max Planck Research School (IMPRS) for Advanced Methods in Process 
and Systems Engineering and Max Planck Society for the Advancement of Science (MPG).
\end{acknowledgments}

\bibliography{Kweyu_refer_aip,mor}

\end{document}